\documentclass{amsart}
\usepackage[utf8]{inputenc}
\usepackage{amsmath}
\usepackage{amssymb}
\usepackage{graphicx}
\usepackage[]{hyperref}
\newtheorem{theorem}{Theorem}[section]
\newtheorem{lemma}[theorem]{Lemma}
\newtheorem{proposition}[theorem]{Proposition}

\theoremstyle{definition}

\newtheorem{remark}[theorem]{Remark}

\numberwithin{equation}{section}

\begin{document}
\title[The quantization coefficient for self-affine measures]
{Convergence order of the quantization error for self-affine measures on Lalley-Gatzouras carpets}

\author{Sanguo Zhu}
\address{School of Mathematics and Physics, Jiangsu University
of Technology\\ Changzhou 213001, China.}
\email{sgzhu@jsut.edu.cn}


\begin{abstract}
 Let $E$ be a Lalley-Gatzouras carpet determined by a set of contractive affine mappings $\{f_{ij}\}_{(i,j)\in G}$. We study the asymptotics of quantization error for the self-affine measures $\mu$ on $E$. We prove that the upper and lower quantization coefficient for $\mu$ are both bounded away from zero and infinity in the exact quantization dimension. This significantly generalizes the previous work concerning the quantization for self-affine measures on Bedford-McMullen carpets. The new ingredients lie in the method to bound the quantization error for $\mu$ from below and that to construct auxiliary measures by applying Prohorov's theorem.
\end{abstract}

\maketitle

\section{Introduction}
The quantization problem consists in the approximation of a given probability measure by discrete probability measures of finite support in $L_r$-metrics.
We refer to Graf and Luschgy \cite{GL:00} for rigorous mathematical foundations of quantization theory and \cite{GN:98} for its deep background in information theory and some engineering technology. One may also see \cite{GL:02,GL:04,GL:05,GL:12,KZ:15,LM:02,PK:01,PK:04,Zhu:23,Zhu:25} for more results on the quantization for fractal measures.
\subsection{Quantization error and its asymptotics}

Let $|x|$ denote the Euclidean norm of $x\in \mathbb{R}^d$. Given $r\in(0,\infty)$, let $\nu$ be a Borel probability measure on $\mathbb{R}^d$ with $\int|x|^rd\mu(x)<\infty$. We denote by ${\rm card}(B)$ the cardinality of a set $B$. For every $n\in\mathbb{N}$, we write $\mathcal{D}_n:=\{\alpha\subset\mathbb{R}^d:1\leq {\rm card}(\alpha)\leq n\}$.
For $x\in\mathbb{R}^d$ and a set $A\subset\mathbb{R}^d$, let $d_A(x):=\inf_{y\in A}d(x,y)$. For a Borel measurable function $g$ on $\mathbb{R}^d$, let $\|g\|_r:=(\int |g|^rd\nu)^{\frac{1}{r}}$. The $n$th quantization error for $\nu$ of order $r$ can be defined by
\[
e_{n,r}(\nu):=\inf_{\alpha\in\mathcal{D}_n}\|d_\alpha(x)\|_r.
\]
By \cite[Lemma 3.4]{GL:00}, $e_{n,r}(\nu)$ is equal to the minimum error in the approximation of $\nu$ with discrete probability measures supported at most $n$ points, in the $L_r$-metric. We refer to \cite{GL:00} for various interpretations of $e_{n,r}(\nu)$ in different contexts.

One of the main goals in quantization theory is to study the asymptotic property of the quantization error, which can be characterized by the $s$-dimensional ($s>0$) upper and lower quantization coefficient:
\[
\overline{Q}_r^s(\nu):=\limsup_{n\to\infty}n^{\frac{r}{s}}e_{n,r}^r(\nu),\;\underline{Q}_r^s(\nu):=\liminf_{n\to\infty}n^{\frac{r}{s}}e_{n,r}^r(\nu).
\]
The upper (lower) quantization dimension $\overline{D}_r(\nu)\; (\underline{D}_r(\nu))$ for $\nu$ of order $r$, is (typically) the critical point at which the upper (lower)
quantization coefficient jumps from infinity to zero (cf. \cite{GL:00,PK:01}). By \cite{GL:00}, we have
\[
\overline{D}_r(\nu)=\limsup_{n\to\infty}\frac{\log n}{-\log e_{n,r}(\nu)};\;\;\underline{D}_r(\nu)=\liminf_{n\to\infty}\frac{\log n}{-\log e_{n,r}(\nu)}.
\]
If $\underline{D}_r(\nu)=\overline{D}_r(\nu)$, we say that the quantization dimension exists and denote the common value by $D_r(\nu)$.

Compared with the upper and lower quantization dimension, we are more concerned about the upper and lower quantization coefficient, because they provide us with exact asymptotic order for the quantization error, when they are both positive and finite.
\subsection{Some known results on the quantization problem}

Let $(f_i)_{i=1}^N$ be a set of contractive similarities on $\mathbb{R}^d$ with contraction ratios $(c_i)_{i=1}^N$. We say that $(f_i)_{i=1}^N$ satisfies the open set condition (OSC), if there exists a bounded non-empty open
set $U$ such that $f_i(U),1\leq i\leq N$, are pairwise disjoint subsets of $U$. According to \cite{Hut:81}, there exists a unique non-empty compact set $F$ satisfying $F=\bigcup_{i=1}^Nf_i(F)$.  We call $F$ the self-similar set determined by $(f_i)_{i=1}^N$. Given a positive probability vector $(p_i)_{i=1}^N$, there exists a unique Borel probability measure satisfying $\nu=\sum_{i=1}^Np_i\nu\circ f_i^{-1}$. This measure is called the self-similar measure associated with $(f_i)_{i=1}^N$ and $(p_i)_{i=1}^N$. Let $\xi_r$ be implicitly defined by $\sum_{i=1}^N(p_ic_i^r)^{\frac{\xi_r}{\xi_r+r}}=1$. Assuming the OSC for $(f_i)_{i=1}^N$, Graf and Luschgy (cf. \cite{GL:02}) proved that
 \[
 0<\underline{Q}^{\xi_r}_r(\nu)\leq\overline{Q}^{\xi_r}_r(\nu)<\infty.
 \]
Moreover, $D_r(\nu)=\xi_r$ increases to the box-counting dimension of $F$ as $r\to\infty$, and decreases to the Hausdorff dimension of $\nu$ as $r\to0$  (cf. \cite{GL:04}).

For every $n\geq 1$, let $\mathcal{C}_n$ denote the partition of $\mathbb{R}^d$ by cubes of the form $\prod_{h=1}^d[k_h2^{-n},(k_h+1)2^{-n})$ with $(k_h)_{h=1}^d\in\mathbb{Z}^d$. For a Borel measure $\nu$, we denote its topological support by $K_\nu$. We define
\[
\mathcal{C}_n^\flat(\nu):=\{C\in\mathcal{C}_n:C\cap K_\nu\neq\emptyset\}.
\]
For every $q\geq 0$, the $L^q$-spectrum $T_\nu(q)$ for $\nu$ can be defined by (cf. \cite{Feng:05,Peres:00})
\[
T_\nu(q):=\lim_{n\to\infty}\frac{\log\sum_{C\in\mathcal{C}_n^\flat(\nu)}\nu(C)^q}{n\log2},
\]
if the limit exists; otherwise, one can consider the upper and lower limits. In the following, we simply write $T(q)$ for $T_\nu(q)$, because no confusion could arise. For $q\neq 1$, the R\'{e}nyi dimension for $\nu$ is defined by $R_\nu(q):=T(q)/(1-q)$.

Kesseb\"ohmer et al identified the upper quantization dimension for an arbitrary compactly supported measure with its R\'{e}nyi dimension at a critical point and provided several sufficient conditions that guarantee the existence of the quantization dimension (cf. \cite[Theorem 1.1]{KNZ:22}). This work, along with Feng-Wang's results in \cite{Feng:05}, yields that the quantization dimensions exist for the self-affine measures on a large class of planar self-affine sets, including Lalley-Gatzouras carpets.
However, the results in \cite{KNZ:22} do not provide us with useful information on the asymptotic order for the quantization error.
In the present paper, we will determine the exact convergence order of the quantization error for the self-affine measures on Lalley-Gatzouras carpets in general.

\subsection{Lalley-Gatzouras carpets}

Let $m\geq 1$ be an integer. For every $1\leq j\leq m$, let $n_j$ be a positive integer. Let $b_j,d_j,1\leq j\leq m$, be positive numbers satisfying
\begin{enumerate}
\item[(A1)] $\sum_{j=1}^mb_j\leq 1;\;d_m\leq 1-b_m$, and

\item[(A2)] $ b_j+d_j\leq d_{j+1}$ for all $1\leq j\leq m-1$ if $m\geq 2$.
\end{enumerate}
For $1\leq j\leq m$, let $a_{ij},1\leq i\leq n_j$, be positive numbers satisfying
\begin{enumerate}
\item[(A3)] $\sum_{i=1}^{n_j}a_{ij}\leq 1,\;\max\limits_{1\leq i\leq n_j}a_{ij}< b_j,a_{n_jj}\leq 1-c_{n_jj}$, and
\item[(A4)] $a_{ij}+c_{ij}\leq c_{(i+1)j}\;\;{\rm for}\;1\leq i\leq n_j-1$ if $n_j\geq 2$.
\end{enumerate}
 Let $G:=\{(i,j):\;1\leq i\leq n_j,1\leq j\leq m\}$ and $G_y:=\{1,\ldots,m\}$. We consider the following mappings on $\mathbb{R}^2$:
\begin{eqnarray}\label{fij}
f_{ij}(x,y)=(x,y)\left(\begin{array}{cccc}
a_{ij} & 0\\
0 &   b_j \\
\end{array}\right)+(c_{ij},d_j),\;(i,j)\in G.
\end{eqnarray}

With the assumptions (A1)-(A4), the iterated function system $\{f_{ij}\}_{(i,j)\in G}$ satisfies the OSC with respect to $U:=(0,1)^2$. By \cite{Hut:81}, there exists a unique non-empty compact set $E$ satisfying
\[
E=\bigcup_{(i,j)\in G}f_{ij}(E).
\]
The set $E$ is called the self-affine set determined by $(f_{ij})_{(i,j)\in G}=:\mathcal{I}$. This type of fractals were first introduced and well studied
by Lalley and Gatzouras \cite{LG:92}, as generalizations of Bedford-McMullen carpets (cf. \cite{Bed:84,Mcmullen:84}. The Hausdorff and box-counting dimension for $E$ were determined; necessary and sufficient conditions were given to guarantee the finiteness and positivity of the Hausdorff measure. One may see \cite{LG:92} for more details.

Given a positive probability vector $\mathcal{P}=(p_{ij})_{(i,j)\in G}$, there exists a unique Borel probability measure $\mu$ satisfying
\begin{equation}\label{prob}
\mu=\sum_{(i,j)\in G}p_{ij}\mu\circ f_{ij}^{-1}.
\end{equation}
We call $\mu$ the self-affine measure associated with $\mathcal{I}$ and $\mathcal{P}$. Let $E_0:=[0,1]^2$. For $l\geq 1$ and $\sigma=((i_1,j_1),\ldots,(i_l,j_l))\in G^l$, we write $f_\sigma:=f_{i_1j_1}\circ\cdots\circ f_{i_lj_l}$. We call $f_\sigma(E_0)$ a \emph{cylinder of order $l$}.

In the past decades, self-affine sets and self-affine measures have attracted great interest of mathematicians
(cf. \cite{Bed:84,Feng:05,KZ:15,King:95,Ko:22,LG:92,Mcmullen:84,Olsen:98,Olsen:11,Peres:94b}). As was noted in the literature, problems concerning typical self-affine sets and self-affine measures are usually difficult. In \cite{Feng:05}, D.-J. Feng and Y. Wang determined various dimensions for a class of planar self-affine sets and self-affine measures, and established several computable formulas for such measures. We refer to \cite{Ko:22,Olsen:11} for related results on self-affine measures on sponges in $\mathbb{R}^d$.

So far, the quantization errors for self-affine measures on Bedford-McMullen carpets have been well studied (cf. \cite{KZ:15,Zhu:18}).
Let $[x]$ denote the largest integer not exceeding $x$ and let $n_0,m_0\in\mathbb{N}$ with $n_0\geq m_0$. Assume that
\begin{eqnarray*}
&&\widetilde{G}\subset\{0,1,\ldots,n_0-1\}\times\{0,1,\ldots,m_0-1\};\\
&&a_{ij}\equiv \frac{1}{n_0}\;((i,j)\in \widetilde{G}),\;b_j\equiv \frac{1}{m_0}\;(j\in \widetilde{G}_y),\\&&c_{ij}\equiv \frac{i}{n_0}\; ((i,j)\in \widetilde{G}),\; d_j\equiv\frac{j}{m_0}\;(j\in \widetilde{G}_y);
\end{eqnarray*}
where $\widetilde{G}_y$ denotes the projection of $\widetilde{G}$ onto the $y$-axis.
Then the carpet $E$ degenerates to a Bedford-McMullen carpet. Let $\mu$ be as defined in (\ref{prob}) with $\widetilde{G}$ in place of $G$. Let $\xi=[\frac{\log m_0}{\log n_0}]$ and $q_j:=\sum_{i:(i,j)\in \widetilde{G}}p_{ij}$. Kesseb\"{o}hmer and Zhu proved that $D_r(\mu)=d_r$, where $d_r$ is given by
\begin{eqnarray}\label{temp1}
\bigg(\sum_{(i,j)\in \widetilde{G}}(p_{ij}m_0^{-r})^{\frac{d_r}{d_r+r}}\bigg)^\xi\bigg(\sum_{j\in \widetilde{G}_y}(q_jm_0^{-r})^{\frac{d_r}{d_r+r}}\bigg)^{1-\xi}=1.
\end{eqnarray}
In \cite{Zhu:18}, Zhu further proved that $0<\underline{Q}^{d_r}_r(\mu)\leq\overline{Q}^{d_r}_r(\mu)<\infty$ holds, in general.

Thanks to the relatively fine structure of Bedford-McMullen carpets, we considered the auxiliary coding space $\widetilde{\Phi}_\infty:=\widetilde{G}^{\mathbb{N}}\times \widetilde{G}_y^{\mathbb{N}}$ and constructed a Bernoulli product $W$ as an auxiliary measure. Then the upper and lower quantization coefficient for $\mu$ were well estimated via the measure $W$, by going back and forth between $E$ and $\Phi_\infty$. In the Lalley-Gatzouras case, however, neither the linear parts nor the translations of the mappings $f_{ij},(i,j)\in G$, need to be constant. This substantially adds to the complexity of the structure of the carpets and the difficulty in analyzing the quantization error.

Combining \cite[Theorem 2]{Feng:05} and \cite[Theorem 1.11]{KZ:15}, we obtain that, the quantization dimension exists for an arbitrary self-affine measure on Lalley-Gatzouras carpets; and the quantization dimension can be expressed in terms of the $L^q$-spectrum of the projection of $\mu$ onto the $y$-axis. However, such expressions seem to be inconvenient in the study of the asymptotic order of the quantization error for $\mu$. We will establish a new expression for the quantization dimension in terms of \emph{cylinders of the same order}. This expression is more convenient for us to construct suitable auxiliary measures which will be used to estimate the quantization coefficient.
\subsection{Some notations}
For $\sigma=((i_1,j_1),\ldots,(i_k,j_k))\in G^k$ and $1\leq h\leq k$, we define
\[
|\sigma|=k,\;\sigma|_h:=((i_1,j_1),\ldots,(i_h,j_h)).
\]
Let $k\geq 1$. Let $\theta$ denote the empty word. For $\sigma\in G^1$, we define $\sigma^-=\theta$. For every $k\geq 2$ and $\sigma\in G^k$, we define $\sigma^-:=\sigma|_{k-1}$.
Let $G^*:=\bigcup_{k\geq 1}G^k$. If $\sigma,\omega\in G^*\cup G^{\mathbb{N}}$ satisfy $\sigma=\omega|_{|\sigma|}$, we say that $\sigma$ is comparable with $\omega$ and write $\sigma\preceq\omega$. If $\sigma\preceq\omega$ and $\sigma\neq\omega$, we write $\sigma\precneqq\omega$.
We call $\sigma,\omega\in G^*$ incomparable if neither $\sigma\preceq\omega$ nor $\omega\preceq\sigma$. For the words $\tau$ in $G_y^*:=\bigcup_{k\geq 1}G_y^k$, we define $\tau|_h,\tau^-$ and the partial order $\preceq$, in a similar manner.

Let $\omega=((i_1,j_1),\ldots,(i_l,j_l))\in G^l$ and $\tau=(j_{l+1},\ldots,j_k)\in G_y^{k-l}$. We write
\begin{eqnarray*}
&&\omega\ast\tau:=((i_1,j_1),\ldots,(i_l,j_l),j_{l+1},\ldots,j_k);\\
&&\omega_y:=(j_1,\ldots,j_l),\;(\omega\ast\tau)_y=(j_1,\ldots,j_k);\\
&&a_\omega:=\prod_{h=1}^la_{i_hj_h},\;b_{\tau}=\prod_{h=l+1}^kb_{j_h};\;a_\theta=b_\theta:=1.
\end{eqnarray*}
For $\sigma=\omega\ast\tau$ with $\omega\in G^*$ and $\tau\in G_y^*$, we write $\omega=:\sigma_L,\tau=:\sigma_R$. Define
\begin{eqnarray*}
\Psi_l:=\big\{\sigma=\sigma_L\ast\sigma_R:b_{\sigma_y^-}\geq a_{\sigma_L}>b_{\sigma_y},\sigma_L\in G^l, \sigma_R\in G_y^*\big\},\;l\geq 1.
\end{eqnarray*}
Let $\Psi^*:=\bigcup_{l\geq 1}\Psi_l$. Let $\sigma=((i_1,j_1),\ldots,(i_l,j_l),j_{l+1},\ldots,j_k)\in\Psi_l$. We define
\begin{eqnarray*}
&&A_{L,\sigma}:=c_{i_1j_1}+\sum_{p=2}^l\bigg(\prod_{h=1}^{p-1}a_{i_hj_h}\bigg)c_{i_pj_p},\;A_{R,\sigma}:=A_{L,\sigma}+\prod_{h=1}^la_{i_hj_h};
\\&&B_{L,\sigma}:=d_{j_1}+\sum_{p=2}^k\bigg(\prod_{h=1}^{p-1}b_{j_h}\bigg)d_{j_p},\;B_{R,\sigma}:=B_{L,\sigma}+\prod_{h=1}^kb_{j_h}.
\end{eqnarray*}
Then a word $\sigma\in \Psi_l$ corresponds to an \emph{approximate square} \emph{$F_\sigma$ }\emph{of order $l$}:
 \[
 F_\sigma:=\big[A_{L,\sigma},A_{R,\sigma}\big]\times\big[B_{L,\sigma},B_{R,\sigma}\big].
 \]
 For $\sigma,\omega\in\Psi^*$, we write $\sigma\preceq\omega$ if $F_\sigma\supset F_\omega$. In Lalley-Gatzouras case, the diameters of approximate squares of the same order can be widely different.

 One can see that $\Psi^*$ provides location codes for approximate squares. However, the way how an approximate square $F_\sigma$ has its sub-approximate squares heavily depends on $\sigma$. This leads to great difficulty in constructing auxiliary measures on the carpets $E$. Therefore, we will also consider the auxiliary coding space $\Phi_\infty=G^\mathbb{N}\times G_y^\mathbb{N}$. Although $\Phi_\infty$ is not completely compatible with the carpet $E$ (cf. Remark \ref{r3}), it suits our purpose perfectly. We define
\begin{eqnarray*}
\Phi_l:=\{\mathcal{L}(\sigma):=\sigma_L\times\sigma_R:\sigma\in\Psi_l\},\;\Phi^*:=\bigcup_{l\geq 1}\Phi_l.
\end{eqnarray*}
We write $\mathcal{L}(\sigma)_L:=\sigma_L$ and $\mathcal{L}(\sigma)_R:=\sigma_R$. For $\sigma,\omega\in\Phi^*$, we write $\sigma\preceq\omega$, if $\sigma_L\preceq\omega_L$ and $\sigma_R\preceq\omega_R$. We call $\sigma,\omega\in\Phi^*$ incomparable if neither $\sigma\preceq\omega$ nor $\omega\preceq\sigma$. For every $\sigma\in\Phi^*$, we define $[\sigma]:=[\sigma_L]\times[\sigma_R]$, where
\[
[\sigma_L]:=\{\rho\in G^{\mathbb{N}}:\sigma_L\preceq\rho\},\;[\sigma_R]:=\{\tau\in G_y^{\mathbb{N}}:\sigma_R\preceq\tau\}.
\]

The main difference between $\Phi^*$ and $\Phi^*$ lies in the fact that they are endowed with different partial orders. Some more information will be presented in Remark \ref{r3}.

For $\sigma,\omega\in\Phi^*$ with $\sigma\preceq\omega$, we have $[\omega]\subset[\sigma]$.
We define
\[
|\sigma_L\ast\sigma_R|=|\sigma_L\times\sigma_R|:=|\sigma_L|+|\sigma_R|,\;(\sigma_L\times\sigma_R)_y:=(\sigma_L\ast\sigma_R)_y.
\]
Let $(p_{ij})_{(i,j)\in G}$ be the same as in (\ref{prob}). For $j\in G_y$, we define $q_j:=\sum_{i=1}^{n_j}p_{ij}$. For
$\sigma_L=((i_1,j_1),\ldots,(i_l,j_l))$ and $\sigma_R:=(j_{l+1},\ldots,j_k)$, we define
\[
p_{\sigma_L}:=\prod_{h=1}^lp_{i_hj_h};\;q_{\sigma_R}:=\prod_{h=l+1}^kq_{j_h}.
\]
\begin{remark}\label{rema}{\rm Let $|B|$ denote the diameter of a set $B\subset\mathbb{R}^2$ and $B^\circ$ its interior. We have the following facts.
\begin{enumerate}
\item[(1)] For every $\sigma\in\Psi^*$, we have $a_{\sigma_L}<|F_\sigma|<\sqrt{2}a_{\sigma_L}$.

\item[(2)] For $\sigma,\omega\in\Psi^*$, we have, either $F_{\sigma}^\circ\cap F_\omega^\circ=\emptyset$, or $F_{\sigma}\subset F_\omega$, or $F_{\omega}\subset F_\sigma$.

\item[(3)] For every $\sigma\in\Psi^*$, it is well known that $\mu(F_\sigma)=p_{\sigma_L}q_{\sigma_R}$.

\item[(4)] The approximate squares as defined above are slightly different from those in \cite{LG:92,Mcmullen:84}. Our definition will enable us to cross out the possibility of the awkward situation that $\sigma_L\precneqq\omega_L$ but $\omega_R\precneqq\sigma_R$. (cf. Lemma \ref{t1}).
\end{enumerate}
}\end{remark}
\subsection{Statement of the main results}

Unlike Bedford-McMullen carpets, we have two major difficulties which arise from the geometric structure of general Lalley-Gatzouras carpets.
\begin{enumerate}
\item[(i)] It is possible that, for some $\sigma,\omega\in\Psi_l$, $|\sigma|\neq|\omega|$, even if $\sigma_L=\omega_L$; for this reason, one can not split the sum $\sum_{\sigma\in\Phi_l}(p_{\sigma_L}q_{\sigma_L}a_{\sigma_L}^r)^t$ into two independent factors, one of which is a sum over $(i,j)\in G$ and the other is a sum over $j\in G_y$ (cf. (\ref{temp1})). This makes the auxiliary measure in \cite{Zhu:18}---a Bernoulli product, fail to work. We will construct auxiliary probability measures by applying Prohorov's theorem.

\item[(ii)] It may happen that, for some $\sigma,\omega\in\Psi^*$ with $\sigma_y\precneqq\omega_y$, we have $|\sigma_L|>|\omega_L|$. This makes the method in
\cite[Lemma 2]{KZ:15} no longer applicable, in case $\min_{j\in G_y}n_j=1$ and $\max_{j\in G_y}n_j\geq 2$. We will present a new method to construct pairwise disjoint approximate squares so that the quantization error for $\mu$ can be estimated from below.
\end{enumerate}
As the main result of the present paper, we will prove
\begin{theorem}\label{mthm}
Let $\mathcal{I}=(f_{ij})_{(i,j)\in G}$ be as defined in (\ref{fij}) and $\mu$ the self-affine measure associated with $\mathcal{I}$ and a positive probability
vector $(p_{ij})_{(i,j)\in G}$. Let $s_r$ be the unique positive number satisfying
\begin{eqnarray}\label{sr}
\lim_{l\to\infty}\frac{1}{l}\log\sum_{\sigma\in\Psi_l}(p_{\sigma_L}q_{\sigma_R}a_{\sigma_L}^r)^{\frac{s_r}{s_r+r}}=0.
\end{eqnarray}
Then we have $0<\underline{Q}_r^{s_r}(\mu)\leq \overline{Q}_r^{s_r}(\mu)<\infty$.
\end{theorem}

\begin{remark}{\rm
(1) The existence of the pressure-like function in (\ref{sr}) and that of the unique number $s_r$ will be proved in Section 3, by applying a sub-additivity argument as used in thermodynamic formulism (cf. \cite{Fal:97}). (2) As a consequence of Theorem \ref{mthm}, we obtain that, $D_r(\mu)=s_r$ and the $n$th quantization error for $\mu$ is of the same order as $n^{-\frac{1}{s_r}}$. This substantially generalizes our previous work on the quantization for self-affine measures on Bedford-McMullen carpets.
}\end{remark}

It was mentioned in \cite[Example 1.15]{KNZ:22}, that for every $q\in (0,1)$, one can obtain an expression for the $L^q$-spectrum $T(q)$ for $\mu$ by means of the formula for the quantization dimension. As an application of Feng-Wang's result \cite[Theorem 2]{Feng:05}, we will show that, for every $q\in[0,\infty)$, $T(q)$ can be expressed in terms of approximate squares of the same order. That is,

\begin{proposition}\label{mthm2}
Let $\mu$ be the same as in Theorem \ref{mthm}. For $q\in[0,\infty)$, the $L^q$-spectrum $T(q)$ for $\mu$ is the unique solution of the following equation:
\[
\lim_{l\to\infty}\frac{1}{l}\log\sum_{\sigma\in\Psi_l}(p_{\sigma_L}q_{\sigma_R})^qa_{\sigma_L}^{T(q)}=0.
\]
\end{proposition}

We will complete the proof of Theorem \ref{mthm} in three steps. First, we will prove some basic facts in Section 2 on the approximate squares. Using these facts, we will construct a set of  pairwise disjoint approximate squares $F_\sigma$ of similar "energy" $\mathcal{E}_r(\sigma):=p_{\sigma_L}q_{\sigma_R}a_{\sigma_L}^r$; then we apply the three-step procedure as described in \cite{KZ:15a} to establish some useful characterizations for the quantization error. Secondly, we will construct some auxiliary measures by applying Prohorov's theorem in section 3. Finally, in section 4, we establish estimates for the quantization coefficient via the auxiliary measures and complete the proof for the main results.
\section{Characterizations of the quantization error}

In the following, we always assume that $m\geq 2$, to avoid trivial cases. For two variables $X,Y$ taking values in $(0,\infty)$,
we write $X\lesssim Y$ ($X\gtrsim Y$), if there exists some constant $C>0$, such that  $X\leq C Y$ ($X\geq CY$).  We write $X\asymp Y$, if we have both $X\lesssim Y$ and $X\gtrsim Y$.

\subsection{Some basic facts}

One can easily construct examples and see that, when $\sigma_L,\omega_L$ are incomparable and $\sigma_y\precneqq\omega_y$, it can happen that $|\sigma_L|>|\sigma_\omega|$. In contrast,
for two words $\sigma,\omega\in\Psi^*$, with $\sigma_L,\omega_L$ comparable, we have
\begin{lemma}\label{tem1}
Let $\sigma,\omega\in\Psi^*$ with $\sigma_L\preceq\omega_L$. Assume that $\sigma_y,\omega_y$ are comparable. Then we have $\sigma_y\preceq\omega_y$; in other words, $|\sigma|\leq|\omega|$.
\end{lemma}
\begin{proof}
Let $\sigma=((i_1,j_1),\ldots,(i_l,j_l),j_{l+1},\ldots,j_k)$. Suppose that $|\sigma|>|\omega|$. Then for some $p_1\geq 0$ and $p_2\geq 1$, we have
\[
\omega=((i_1,j_1),\ldots,(i_l,j_l),\ldots,(i_{l+p_1},j_{l+p_1}),j_{l+p_1+1},\ldots,j_{k-p_2}).
\]
By the definition of $\Psi^*$, we have, $\prod_{h=1}^{k-1}b_{j_h}\geq\prod_{h=1}^la_{i_hj_h}>\prod_{h=1}^kb_{j_h}$. Hence,
\[
\prod_{h=1}^{k-p_2}b_{j_h}\geq\prod_{h=1}^{k-1}b_{j_h}\geq\prod_{h=1}^la_{i_hj_h}\geq\prod_{h=1}^{l+p_1}a_{i_hj_h}.
\]
This contradicts the fact that $\omega\in\Psi_{l+p_1}$.
\end{proof}

For convenience, we define
\begin{eqnarray}\label{temp2}
&\underline{a}=\min\limits_{(i,j)\in G}a_{ij},\;\overline{a}=\max\limits_{(i,j)\in G}a_{ij};\;\; \underline{b}=\min\limits_{j\in G_y}b_j,\;\overline{b}=\max\limits_{j\in G_y}b_j;\nonumber\\
&A_1=[\frac{\log\underline{a}}{\log \overline{b}}]+1,\;A_2=[\frac{\log\underline{b}}{\log \overline{a}}]+1,\;A_3=\max\limits_{j\in G_y}\max\limits_{1\leq i\leq n_j}\frac{a_{ij}}{b_j},\;A_4=\frac{\log A_3}{\log\underline{b}}.
\end{eqnarray}
These constants will be used to control the length of $\tau_R$ in the construction of sub-approximate squares $F_{\tau_L\ast\tau_R}$.
\begin{lemma}\label{s4}
Assume that $\sigma,\omega\in\Psi^*$ and $F_\omega\subseteq F_\sigma$. Then we have
\begin{enumerate}
\item[(i)] if $|\omega_L|=|\sigma_L|+1$, we have $0\leq|\omega|-|\sigma|\leq A_1$;
\item[(ii)] if $|\omega_L|-|\sigma_L|\geq A_2$, then we have, $|\omega|\geq|\sigma|+1$.
\end{enumerate}
\end{lemma}
\begin{proof}
We write $\sigma=((i_1,j_1),\ldots,(i_l,j_l),j_{l+1},\ldots,j_k)$. Since $F_\omega\subseteq F_\sigma$, we have, $\sigma_L\preceq\omega_L$ and $\sigma_y\preceq\omega_y$. It follows that $|\omega|\geq|\sigma|$.

(i) Assume that $|\omega_L|=|\sigma_L|+1$. Then for some $1\leq i\leq n_{j_{l+1}}$, we have,
\[
\omega_L=((i_1,j_1),\ldots,(i_l,j_l),(i,j_{l+1}),j_{l+2},\ldots, j_k,\ldots, j_{k+h}).
\]
Suppose that $h>A_1(\geq 2)$. By the definition of $\Psi_l$, we deduce
\[
\prod_{p=1}^{k+h-1}b_{j_p}=\prod_{p=1}^kb_{j_p}\prod_{p=k+1}^{k+h-1}b_{j_p}\leq \prod_{p=1}^la_{i_pj_p}\overline{b}^{A_1}
<\prod_{p=1}^la_{i_pj_p}\underline{a}\leq\prod_{p=1}^{l+1}a_{i_pj_p}.
\]
This contradicts the fact that $\omega\in\Psi_{l+1}$. Hence, $h\leq A_1$ and $|\omega|\leq k+A_1$.

(ii) Assume that $\omega_L=((i_1,j_1),\ldots,(i_l,j_l),(i,j_{l+1}),\ldots,(i_{l+h},j_{l+h}))$ with $h\geq A_2$. Again, by the definition of $\Psi_l$, we have
\[
a_{\sigma_L}\prod_{p=l+1}^{l+h}a_{i_pj_p}\leq a_{\sigma_L}\overline{a}^{A_2}<\prod_{p=1}^{k-1}b_{j_p}\underline{b}\leq\prod_{p=1}^kb_{j_p}.
\]
It follows that $|\omega|\geq|\sigma|+1$.
\end{proof}
Let $\sigma=((i_1,j_1),\ldots,(i_l,j_l),j_{l+1},\ldots,j_k)\in\Psi_l$. For $\omega\in\Psi_{l+1}$ with $\sigma\preceq\omega$, by Lemma \ref{s4}, for some $1\leq i\leq n_{j_{l+1}}$ and $p\geq 0$, we have
\[
\omega=((i_1,j_1),\ldots,(i_l,j_l),(i,j_{l+1}),j_{l+2},\ldots,j_k,\ldots,j_{k+p}).
\]

 We end this subsection with the following observation on the length of $\sigma_R$ for a word $\sigma\in\Psi^*$. We will need it in the characterization for $e_{\varphi_{n,r},r}(\mu)$.
\begin{lemma}\label{zzz1}
For $\sigma\in\Psi_1$, we have $|\sigma_R|\geq 1$; for $\sigma\in\Psi_l$, we have $|\sigma_R|\geq A_4l$.
\end{lemma}
\begin{proof}
For $(i,j)\in G$, we have $b_j>a_{ij}$. This implies that $|\sigma_R|\geq 1$ for every $\sigma\in\Psi_1$.
Let $\sigma=((i_1,j_1),\ldots,(i_l,j_l),j_{l+1},\ldots,j_k)\in\Psi_l$. We have
\begin{eqnarray*}
\underline{b}^{k-l}\leq\prod_{h=l+1}^kb_{j_h}=\prod_{h=1}^kb_{j_h}\prod_{h=1}^lb_{j_h}^{-1}<\prod_{h=1}^l\frac{a_{i_hj_h}}{b_{j_h}}\leq A_3^l.
\end{eqnarray*}
It follows that $|\sigma_R|=k-l\geq A_4l$.
\end{proof}

\subsection{Characterizations for the quantization error}
We call a finite subset $\Gamma$ of $\Psi^*$ a finite anti-chain, if $F_\sigma, \sigma\in\Gamma$, are pairwise non-overlapping; if, in addition,
$E\subset\bigcup_{\sigma\in\Gamma}F_\sigma$, then we call $\Gamma$ a finite maximal anti-chain in $\Psi^*$. For $A,B\subset\mathbb{R}^2$, we define
\[
d_h(A,B):=\inf_{(x,y)\in A,(x',y')\in B}|x-x'|,\;d_v(A,B):=\inf_{(x,y)\in A,(x',y')\in B}|y-y'|.
\]

In order to establish estimates for the quantization error, we will construct a finite anti-chain $\mathcal{F}_{n,r}$, such that $\mathcal{E}_r(\tau)\asymp\underline{\eta}_r^n$ for every $\tau\in \mathcal{F}_{n,r}$, and for every pair of distinct words $\tau^{(1)},\tau^{(2)}$ in $\mathcal{F}_{n,r}$, the following holds:
\begin{equation}\label{gg4}
d(F_{\tau^{(1)}},F_{\tau^{(2)}})\gtrsim\max\{|F_{\tau^{(1)}}|,|F_{\tau^{(2)}}|\}.
\end{equation}

For every $\sigma\in\Psi_1$, we define $\sigma^\flat:=\theta$. Next, we assume that $l\geq 2$. For $\sigma\in\Psi_l$, we write $\sigma=((i_1,j_1),\ldots,(i_l,j_l),j_{l+1},\ldots,j_k)$. We have
\[
\prod_{h=1}^{l-1}a_{i_hj_h}>a_{i_lj_l}^{-1}\prod_{h=1}^kb_{j_h}>\prod_{h=1}^kb_{j_h}.
\]
Thus, there exists a unique integer $p\geq 0$, such that
\[
\prod_{h=1}^{k-p-1}b_{j_h}\geq \prod_{h=1}^{l-1}a_{i_hj_h}>\prod_{h=1}^{k-p}b_{j_h}.
\]
We define $\sigma^\flat:=\sigma_L^-\ast (j_l,\ldots,j_{k-p})$. One can see that $F_\sigma\subset F_{\sigma^\flat}$. Write
\begin{eqnarray}\label{energy}
&&\underline{p}:=\min_{(i,j)\in G}p_{ij},\;\overline{p}:=\max_{(i,j)\in G}p_{ij},\;\underline{q}:=\min_{j\in G_y}q_j,\;\overline{q}:=\max_{j\in G_y}q_j;\nonumber\\
&&\mathcal{E}_r(\theta):=1;\;\mathcal{E}_r(\sigma):=\mu(F_\sigma)a_{\sigma_L}^r;\;\overline{\eta}_r:=\overline{a}^r\max_{(i,j)\in G}\frac{p_{ij}}{q_j};\;
\underline{\eta}_r:=\underline{p}\;\underline{q}^{A_1}\underline{a}^{3r}.
\end{eqnarray}
From Remark \ref{rema} (2) and Lemma \ref{s4}, for every $\sigma\in\Psi^*$, we have
\begin{equation}\label{z3}
\underline{\eta}_r\mathcal{E}_r(\sigma^\flat)<\underline{q}^{A_1}\underline{a}^{r}\mathcal{E}_r(\sigma^\flat)\leq\mathcal{E}_r(\sigma)\leq \overline{\eta}_r\mathcal{E}_r(\sigma^\flat)<\mathcal{E}_r(\sigma^\flat).
\end{equation}
For every $n\geq 1$, (\ref{z3}) allows us to define
\begin{eqnarray}\label{lambdajr}
\Lambda_{n,r}:=\{\sigma\in\Psi^*:\mathcal{E}_r(\sigma^\flat)\geq\underline{\eta}_r^n>\mathcal{E}_r(\sigma)\};\;\varphi_{n,r}:={\rm card}(\Lambda_{n,r}).
\end{eqnarray}
\begin{remark}\label{cardlambdajr}{\rm
(1) For every $n\geq 1$, $\Lambda_{n,r}$ is a finite maximal anti-chain in $\Psi^*$. As we did in \cite[Lemma 1]{KZ:15}, it is easy to show that $\varphi_{n,r}\asymp\varphi_{n+1}$. Moreover, for every $\sigma\in\Lambda_{n,r}$, we have, $\underline{\eta}_r^{n+1}\leq \mathcal{E}_r(\sigma)<\underline{\eta}_r^n$. (2) Let $l_{n,r}:=\min_{\sigma\in\Lambda_{n,r}}|\sigma_L|$. We have, $\underline{\eta}_r^{l_{n,r}}<\underline{\eta}_r^n$, so $l_{n,r}\geq n$.
}\end{remark}
\begin{remark}\label{sigmalincomparable}{\rm
\rm For $\sigma,\omega\in\Lambda_{n,r}$, if $\sigma_L,\omega_L$ are comparable, then $\sigma_y,\omega_y$ are incomparable, because otherwise, by Lemma \ref{tem1}, we would have $F_\sigma\subseteq F_\omega$, or $F_\omega\subseteq F_\sigma$, contradicting the fact that $\Lambda_{n,r}$ is a finite anti-chain.
}\end{remark}

Next, we are going to construct the finite anti-chain $\mathcal{F}_{n,r}\subset\Psi^*$ out of $\Lambda_{n,r}$.
First, let us cope with the extreme case that $n_j=1$ for every $j\in G_y$, where, $\mathcal{F}_{n,r}$ can be defined in a convenient manner.
\begin{lemma}\label{zs1}
Assume that $n_j=1$ for every $j\in G_y$. Let $\sigma,\omega\in\Psi^*$ with $\sigma_y,\omega_y$ comparable. Then we have, either $F_\omega\subseteq F_\sigma$, or $F_\sigma\subseteq F_\omega$. In particular, for every pair $\sigma,\omega$ of distinct words in $\Lambda_{n,r}$, $\sigma_y,\omega_y$ are incomparable.
\end{lemma}
\begin{proof}
By the assumptions that $\sigma_y, \omega_y$ are comparable and $n_j=1$ for all $j\in G_y$, we know that $\sigma_L,\omega_L$ are comparable. The lemma follows from Lemma \ref{tem1}.
\end{proof}
\begin{remark}{\rm
Let $n>2A_4^{-1}A_2$. Assume that $n_j=1$ for every $j\in G_y$. For $\sigma\in\Lambda_{n,r}$, we write
$\sigma=((1,j_1),\ldots,(1,j_l),j_{l+1},\ldots,j_k)$. From Remark \ref{cardlambdajr} (2) and Lemma \ref{zzz1}, we have, $|\sigma_R|>2A_2$. By Lemma \ref{s4} (ii), we may define
\[
\tau_\sigma:=\sigma_L\ast((1,j_{l+1}),\ldots,(1,j_{l+2A2}))\ast(j_{l+2A_2+1},\ldots,j_k,1,m,\ldots,j_{\overline{k}})\in\Psi^*.
\]
Then for $\sigma,\omega\in\Lambda_{n,r}$ with $\sigma\neq\omega$, we have
\[ d_v(F_{\tau_\sigma},F_{\tau_\omega})\geq\underline{b}^2\max\{b_{\sigma_y},b_{\omega_y}\}\geq\underline{b}^2(1+\underline{b}^{-2})^{-\frac{1}{2}}
\max\{|F_{\tau_\sigma}|,|F_{\tau_\omega}|\}.
\]
Thus, in case that $n_j=1$ for every $j\in G_y$, it is sufficient to define
\[
\mathcal{F}_{n,r}:=\{\tau_\sigma:\sigma\in\Lambda_{n,r}\}.
\]}
\end{remark}

In the following, we assume that $n_{j_0}\geq 2$ for some $j_0\in G_y$. We will construct $\mathcal{F}_{n,r}$ in two steps. For the choice of $\overline{\sigma}$, we are inspired by \cite[Lemma 2]{KZ:15} and \cite[Lemma 7.2]{Ko:22}.

\textbf{Step 1}: For $\sigma=((i_1,j_1),\ldots,(i_l,j_l),j_{l+1},\ldots,j_k)\in\Lambda_{n,r}$, we define
\begin{equation}\label{gg3}
\overline{\sigma}:=((i_1,j_1),\ldots,(i_l,j_l),(1,j_0),(n_{j_0},j_0),j_{l+1},\ldots,j_k,\ldots, j_{\overline{k}})\in\Psi^*.
\end{equation}
We define $\mathcal{B}_{n,r}:=\{\overline{\sigma}:\sigma\in\Lambda_{n,r}\}$. We will prove that, for all large $n$, for every pair of distinct words $\overline{\sigma},\overline{\omega}\in\mathcal{B}_{n,r}$, we have, either of the following holds:

(i) $d_h(F_{\overline{\sigma}},F_{\overline{\omega}})\gtrsim\max\{|F_{\overline{\sigma}}|,|F_{\overline{\omega}}|\}$ (Lemma \ref{ss4});

(ii) $\overline{\sigma}_y,\overline{\omega}_y$ are incomparable (Lemma \ref{zsg1}).

\textbf{Step 2}: For every $\overline{\sigma}\in \mathcal{B}_{n,r}$, we select a word $\overline{\sigma}^*\in\Psi^*$ such that $F_{\overline{\sigma}^*}\subset F_{\overline{\sigma}}$ and 
\[
d_v(F_{\overline{\sigma}^*},F_{\overline{\omega}^*})\gtrsim\max\{|F_{\overline{\sigma}^*}|,|F_{\overline{\sigma}^*}|\},
 \]
whenever $\overline{\sigma}_y,\overline{\omega}_y$ are incomparable. We define $\mathcal{F}_{n,r}:=\{\overline{\sigma}^*:\overline{\sigma}\in\mathcal{B}_{n,r}\}$.

\begin{lemma}\label{ss4}
Let $\sigma,\omega\in\Lambda_{n,r}$ with $\sigma_L,\omega_L$ incomparable and $\overline{\sigma}_y,\overline{\omega}_y$ comparable. We have,
$d_h(F_{\overline{\sigma}},F_{\overline{\omega}})\geq2^{-\frac{1}{2}}\underline{a}^2\max\{|F_{\overline{\sigma}}|,|F_{\overline{\omega}}|\}$.
\end{lemma}
\begin{proof}
Let $\tau_0=((1,j_0),(n_{j_0},j_0))$. Without loss of generality, we assume that $|\sigma_L|=l_1,|\omega_L|=l_2$ and that $l_1\leq l_2$. Since $\overline{\sigma}_y,\overline{\omega}_y$ are comparable, we have that $(\omega_L)_y|_{l_1}=(\sigma_L)_y$. As $\sigma_L,\omega_L$ are incomparable, we know that $\omega_L|_{l_1}\neq\sigma_L$. Therefore, $f_{\sigma_L}(E_0),f_{\omega_L|_{l_1}}(E_0)$ are non-overlapping. Note that $F_{\overline{\sigma}}\subset f_{\sigma_L\ast\tau_0}(E_0)$ and $F_{\overline{\omega}}\subset f_{\omega_L\ast\tau_0}(E_0)\subset f_{\omega_L|_{l_1}}(E_0)$. We deduce
\begin{eqnarray*}
d_h(F_{\overline{\sigma}},F_{\overline{\omega}})&\geq& d_h(f_{\sigma_L\ast\tau_0}(E_0),f_{\omega_L\ast\tau_0}(E_0))\\&\geq&\underline{a}^2\max\{a_{\sigma_L},a_{\omega_L}\}
\geq2^{-\frac{1}{2}}\underline{a}^2\max\{|F_{\overline{\sigma}}|,|F_{\overline{\omega}}|\}.
\end{eqnarray*}
This completes the proof of the lemma.
\end{proof}

In the following, we are going to show that, for $\sigma,\omega\in\Lambda_{n,r}$ with $\sigma_L\preceq\omega_L$, $\overline{\sigma}_y,\overline{\omega}_y$, are necessarily incomparable. We begin with the simplest cases when $|\omega_L|-|\sigma_L|\leq 2$. By Lemma \ref{zzz1} and Remark \ref{cardlambdajr} (2), for every $n\geq 4A_4^{-1}$ and every $\tau\in\Lambda_{n,r}$, we have $|\tau_R|\geq 4$. We assume that $n\geq 4A_4^{-1}$ in the subsequent Lemma \ref{zz0} in order to avoid some trivial cases.
\begin{lemma}\label{zz0}
Let $\sigma,\omega\in\Lambda_{n,r}$ be distinct words. Assume that $\sigma_L\preceq\omega_L$ and $|\omega_L|-|\sigma_L|\leq 2$, then $\overline{\sigma}_y,\overline{\omega}_y$ are incomparable.
\end{lemma}
\begin{proof}
By Remark \ref{sigmalincomparable}, when $\sigma_L\preceq\omega_L$, we have $\sigma_y,\omega_y$ are incomparable. We assume that $|\sigma|=k_1$ and $|\omega|=k_2$ and write (for $0\leq h\leq 2$)
\begin{eqnarray}\label{sigmaomega1}
&&\sigma=((i_1,j_1),\ldots,(i_l,j_l),j_{l+1},\ldots,j_{k_1});\\
&&\omega=((i_1,j_1),\ldots,(i_l,j_l),\ldots,(\hat{i}_{l+h},\hat{j}_{l+h}),\hat{j}_{l+h+1},\ldots,\hat{j}_{k_2}).\label{sigmaomega2}
\end{eqnarray}

(i) First we assume that $\sigma_L=\omega_L$, then  $\sigma_R,\omega_R$ are incomparable. Note that
$\overline{\sigma}_y|_{k_1+2}=(\sigma_L)_y\ast (j_0,j_0)\ast\sigma_R$ and $\overline{\omega}_y|_{k_2+2}=(\sigma_L)_y\ast (j_0,j_0)\ast\omega_R$.
It follows that $\overline{\sigma}_y,\overline{\omega}_y$ are incomparable.

(ii) Now we assume that $|\omega_L|=|\sigma_L|+1$. We write
\begin{eqnarray*}
&&\overline{\sigma}:=\sigma_L\ast\tau_0\ast(j_{l+1},j_{l+2},\ldots,j_{k_1},\ldots, j_{\overline{k}_1});\\
&&\overline{\omega}:=\sigma_L\ast(\hat{i}_{l+1},\hat{j}_{l+1})\ast\tau_0\ast(\hat{j}_{l+2},\ldots,\hat{j}_{k_2},\ldots, \hat{j}_{\overline{k}_2}).
\end{eqnarray*}
For convenience, we write $\overline{\sigma}_y$ and $\overline{\omega}_y$ in detail:
\begin{eqnarray*}
&&\overline{\sigma}_y=(j_1,\ldots,j_l,j_0,j_0,j_{l+1},j_{l+2},\ldots,j_{k_1},\ldots,j_{\overline{k}_1}),\\
&&\overline{\omega}_y=(j_1,\ldots,j_l,\hat{j}_{l+1},j_0,j_0,\hat{j}_{l+2},\ldots,\hat{j}_{k_2},\ldots,\hat{j}_{\overline{k}_2}).
\end{eqnarray*}
If $\hat{j}_{l+1}\neq j_0$, then $\overline{\sigma}_y,\overline{\omega}_y$ are incomparable. Next, we assume that $\hat{j}_{l+1}=j_0$. If $j_{l+1}\neq j_0$, we again have that $\overline{\sigma}_y,\overline{\omega}_y$ are incomparable; otherwise, we have, $\hat{j}_{l+1}=j_{l+1}=j_0$. Note that $\sigma_y,\omega_y$ are incomparable. Thus, $(j_{l+2},\ldots,j_{k_1})$ and $(\hat{j}_{l+2},\ldots,\hat{j}_{k_2})$ are incomparable. Hence, $\overline{\sigma}_y,\overline{\omega}_y$ are incomparable.

(iii) The case that $|\omega_L|=|\sigma_L|+2$ can be proved analogously to (ii).
\end{proof}

Next, we prove by contradiction that, for $\sigma,\omega\in\Lambda_{n,r}$ with $\sigma_L\precneqq\omega_L$ and $|\omega_L|-|\sigma_L|\geq3$, $\overline{\sigma}_y, \overline{\omega}_y$ are incomparable: assuming the contrary, we will construct some word $\breve{\omega}\precneqq\omega$ such that $\breve{\omega}\in\Lambda_{n,r}$, contradicting the fact that $\Lambda_{n,r}$ is a finite anti-chain.

\begin{lemma}\label{zsg1}
Let $A_{5,r}:=\big[\frac{\log(\underline{q}^2\underline{\eta}_r)}{r\log \overline{a}}\big]$ and $n>T_{1,r}:=[A_4^{-1}(A_1+A_{5,r}+3)]$. Let $\sigma, \omega\in\Lambda_{n,r}$. Assume that $\sigma_L\precneqq\omega_L$ and $|\omega_L|-|\sigma_L|\geq 3$. Then $\overline{\sigma}_y, \overline{\omega}_y$ are incomparable.
\end{lemma}
\begin{proof}
Assume that $\sigma,\omega\in\Lambda_{n,r},\sigma_L\precneqq\omega_L$, but $\overline{\sigma}_y,\overline{\omega}_y$ are comparable. We will deduce a contradiction.
Let $\sigma,\omega$ be the same as in (\ref{sigmaomega1}) and (\ref{sigmaomega2}) with $h\geq 3$. Then we have $\overline{\sigma}:=\sigma_L\ast\tau_0\ast(j_{l+1},\ldots,j_{k_1},\ldots, j_{\overline{k}_1})$ and
\begin{eqnarray*}
\overline{\omega}:=\sigma_L\ast((\hat{i}_{l+1},\hat{j}_{l+1}),\hat{i}_{l+2},\hat{j}_{l+2})\ldots,(\hat{i}_{l+h},\hat{j}_{l+h}))\ast\tau_0
\ast(\hat{j}_{l+h+1},\ldots,\hat{j}_{k_2},\ldots, \hat{j}_{\overline{k}_2}).
\end{eqnarray*}

\textbf{Claim 1}: we have $k_2\geq k_1$. Suppose that $k_2<k_1$. Then $k_1>l+h$ and
\begin{eqnarray*}
&&\overline{\sigma}_y=(j_1,\ldots,j_l,j_0,j_0,j_{l+1},\ldots,j_{l+h-2},j_{l+h-1},j_{l+h},j_{l+h+1}\ldots,j_{k_1},\ldots,j_{\overline{k}_1}),\\
&&\overline{\omega}_y=(j_1,\ldots,j_l,\hat{j}_{l+1},\hat{j}_{l+2},\hat{j}_{l+3},\hat{j}_{l+4},\ldots,\hat{j}_{l+h},j_0,j_0,\hat{j}_{l+h+1},\ldots,\hat{j}_{k_2},\ldots,
\hat{j}_{\overline{k}_2}).
\end{eqnarray*}
By the assumption, we have that $\overline{\sigma}_y,\overline{\omega}_y$ are comparable. Hence,
\begin{eqnarray}
&&(\hat{j}_{l+1},\hat{j}_{l+2})=(j_0,j_0)=(j_{l+h-1},j_{l+h});\label{ggg0}\\&&(j_{l+1},\ldots,j_{l+h-2})=(\hat{j}_{l+3},\ldots,\hat{j}_{l+h});\label{ggg1}\\
&&(\hat{j}_{l+h+1},\ldots,\hat{j}_{k_2})\precneqq(j_{l+h+1},\ldots,j_{k_1}).\nonumber
\end{eqnarray}
Note that $\sigma\in\Psi_l$. we have $\prod_{p=1}^{k_1-1}b_{j_p}\geq a_{\sigma_L}>\prod_{p=1}^{k_1}b_{j_p}$. Hence,
\begin{eqnarray*}
b_{\omega_y}&=&\prod_{p=1}^lb_{j_p}\prod_{p=l+1}^{l+h}b_{\hat{j}_p}\prod_{p=l+h+1}^{k_2}b_{\hat{j}_p}\\
&=&\prod_{p=1}^lb_{j_p}\bigg(b_{j_0}^2\prod_{p=l+3}^{l+h}b_{\hat{j}_p}\bigg)\prod_{p=l+h+1}^{k_2}b_{\hat{j}_p}\;\;({\rm by}\;(\ref{ggg0}))\\
&=&\prod_{p=1}^lb_{j_p}\bigg(\prod_{p=l+1}^{l+h-2}b_{j_p}b_{j_0}^2\bigg)\prod_{p=l+h+1}^{k_2}b_{j_p}\;\;({\rm by}\;(\ref{ggg1}))\\
&\geq&\prod_{p=1}^{k_1-1}b_{j_p}\geq a_{\sigma_L}>a_{\omega_L}.
\end{eqnarray*}
This contradicts the fact that $\omega\in\Psi_{l+h}$ and Claim 1 follows. Thus,
\begin{equation}\label{ggg2}
(j_{l+h+1},\ldots,j_{k_1})\preceq(\hat{j}_{l+h+1},\ldots,\hat{j}_{k_2})\;\;{\rm if}\;k_1>l+h.
\end{equation}

\textbf{Claim 2}: we have $\mathcal{E}_r(\omega)\leq\underline{q}^{-2}\overline{a}^{hr}\mathcal{E}_r(\sigma)$. We distinguish two cases.

Case 1: $k_1>l+h$. In this case, using (\ref{ggg0})-(\ref{ggg2}), we deduce
\begin{eqnarray*}
\mathcal{E}_r(\omega)&=&\big(p_{\sigma_L}\prod_{p=l+1}^{l+h}p_{\hat{i}_p\hat{j}_p}\prod_{p=l+h+1}^{k_2}q_{\hat{j}_p}\big)(a_{\sigma_L}^r
\prod_{p=l+1}^{l+h}a_{\hat{i}_p\hat{j}_p}^r)
\\&\leq&\big(p_{\sigma_L}p_{\hat{i}_{l+1}\hat{j}_{l+1}}p_{\hat{i}_{l+2}\hat{j}_{l+2}}\prod_{p=l+3}^{l+h}q_{\hat{j}_p}\prod_{p=l+h+1}^{k_2}q_{\hat{j}_p}\big)(a_{\sigma_L}^r
\prod_{p=l+1}^{l+h}a_{\hat{i}_p\hat{j}_p}^r)\\
&<&\big(p_{\sigma_L}\prod_{p=l+1}^{l+h-2}q_{j_p}\prod_{p=l+h+1}^{k_1}q_{j_p}\big)(a_{\sigma_L}^r\prod_{p=l+1}^{l+h}a_{\hat{i}_p\hat{j}_p}^r)\;\;({\rm By}\;(\ref{ggg1}))\\
&\leq&\underline{q}^{-2}\overline{a}^{hr}\mathcal{E}_r(\sigma).
\end{eqnarray*}

Case 2: $k_1\leq l+h$. If $k_1\geq l+h-2$, then (\ref{ggg1}) remains valid; otherwise, we have
$(j_{l+1},\ldots,j_{k_1})=(\hat{j}_{l+3},\ldots,\hat{j}_{k_1+2})$. Let $k:=\min(k_1,l+h-2)$, we deduce
\begin{eqnarray*}
\mathcal{E}_r(\omega)&\leq&\big(p_{\sigma_L}\prod_{p=l+1}^{k}q_{j_p}\big)(a_{\sigma_L}^r\prod_{p=l+1}^{l+h}a_{\hat{i}_p\hat{j}_p}^r)\\
&\leq&\big(p_{\sigma_L}\prod_{p=l+1}^{k_1-2}q_{j_p}\big)(a_{\sigma_L}^r\prod_{p=l+1}^{l+h}a_{\hat{i}_p\hat{j}_p}^r)\\
&\leq&\underline{q}^{-2}\overline{a}^{hr}\mathcal{E}_r(\sigma).
\end{eqnarray*}

\textbf{Claim 3}: we have  $h\leq A_{5,r}$.
Assume that $h>A_{5,r}$. Then by Claim 2, we have $\mathcal{E}_r(\omega)<\underline{\eta}_r^{n+1}$, a contradiction (cf. Remark \ref{cardlambdajr}).

Using Claims 1-3, we are able to complete the proof of the lemma. Because $n>T_{1,r}$, by Lemma \ref{zzz1} and Remark \ref{cardlambdajr} (2), we have,
\begin{equation}\label{temp4}
|\sigma_R|=k_1-l\geq A_1+A_{5,r}+3\geq A_1+h+3.
\end{equation}
Thus, we may write $\sigma_y,\omega_y$ in detail:
\begin{eqnarray*}
&&\sigma_y=(j_1,\ldots,j_l,j_{l+1},j_{l+2},\ldots,j_{l+h-2},j_0,j_0,j_{l+h+1}\ldots,j_{k_1}),\\
&&\omega_y=(j_1,\ldots,j_l,j_0,j_0,\hat{j}_{l+3},\hat{j}_{l+4},\ldots,\hat{j}_{l+h},\hat{j}_{l+h+1},\ldots,\hat{j}_{k_2}).
\end{eqnarray*}
Let $\breve{\omega}:=\sigma_L\ast(j_0,j_0,\hat{j}_{l+3},\ldots,\hat{j}_{l+h},\hat{j}_{l+h+1},\ldots,\hat{j}_{k_1})$. By (\ref{ggg0})-(\ref{temp4}), $(\breve{\omega})_R\; ({\rm respectively}\;(\breve{\omega})_R^-)$ is a permutation of $\sigma_R\;({\rm respectively}\;\sigma_R^-)$. Thus,
\begin{eqnarray*}
b_{\breve{\omega}_y^-}=b_{\sigma_y^-}\geq a_{\sigma_L}=a_{\breve{\omega}_L}>b_{\sigma_y}=b_{\breve{\omega}_y}.
\end{eqnarray*}
It follows that $\breve{\omega}\in\Psi_l$ and $F_\omega\subsetneq F_{\breve{\omega}}$. Note that, by Lemma \ref{s4} (ii), we have, $|(\sigma^\flat)_R|\geq |\sigma_R|-A_1\geq A_{5,r}+3\geq h+3$. For some integer $T\geq 2$, we may write
\[
(\sigma^\flat)_R=(j_l,j_{l+1},\ldots,j_{l+h-2},j_0,j_0,j_{l+h+1},\ldots,j_{l+h+T}).
\]
Again, by (\ref{ggg0})-(\ref{temp4}) and the definition of $\Psi_{l-1}$, we deduce
\[
(\breve{\omega}^\flat)_R=(j_l,j_0,j_0,\hat{j}_{l+3},\ldots,\hat{j}_{l+h},\hat{j}_{l+h+1},\ldots,\hat{j}_{l+h+T}).
\]
Thus, $(\breve{\omega}^\flat)_R$ is also a permutation of $(\sigma^\flat)_R$. It follows that
\begin{eqnarray*}
\mathcal{E}_r(\breve{\omega}^\flat)=\mathcal{E}_r(\sigma^\flat)\geq\underline{\eta}_r^n>\mathcal{E}_r(\sigma)=\mathcal{E}_r(\breve{\omega}).
\end{eqnarray*}
This implies that $\breve{\omega}\in\Lambda_{n,r}$ and $\omega\notin\Lambda_{n,r}$, contradicting the hypothesis.
\end{proof}

  In the remaining part of this section, we assume that $n>T_{2,r}:=T_{1,r}+2A_4^{-1}A_2$. Let $\overline{\sigma}\in\mathcal{B}_{n,r}$ be as defined in (\ref{gg3}). For every $1\leq h\leq 2A_2$ , we fix an integer $i_{l+h}\in [1,n_{j_{l+h}}]$ and define (cf. Lemma \ref{s4})
\begin{eqnarray*}
&&\overline{\sigma}^*_L:=\sigma_L\ast\tau_0\ast\big((i_{l+1},j_{l+1}),\ldots,(i_{l+2A_2},j_{l+2A_2})\big),\\
&&\overline{\sigma}^*_R:=(j_{l+2A_2+1},\ldots,j_k,\ldots, j_{\overline{k}},1,m,\ldots,j_{\widetilde{k}}),
\end{eqnarray*}
such that $\overline{\sigma}^*:=\overline{\sigma}^*_L\ast\overline{\sigma}^*_R\in\Psi^*$. Then we have $F_{\overline{\sigma}^*}\subset F_{\overline{\sigma}}$. We define
\begin{equation}\label{Fjr}
\mathcal{F}_{n,r}:=\{\overline{\sigma}^*:\overline{\sigma}\in\mathcal{B}_{n,r}\}.
\end{equation}

\begin{lemma}\label{zz2}
For every pair $\overline{\sigma}^*,\overline{\omega}^*$ of distinct words in $\mathcal{F}_{n,r}$, we have
\[
d(F_{\overline{\sigma}^*},F_{\overline{\omega}^*})\geq(1+\underline{b}^{-2})^{-1}\underline{a}^2\max\{|F_{\overline{\sigma}^*}|,|F_{\overline{\omega}^*}|\}.
\]
\end{lemma}
\begin{proof}
If $(\overline{\sigma})_y,(\overline{\omega})_y$ are incomparable, so are $(\overline{\sigma}^*)_y,(\overline{\omega}^*)_y$. We have
\begin{eqnarray*}
d(F_{\overline{\sigma}^*},F_{\overline{\omega}^*})&\geq& d_v(F_{\overline{\sigma}^*},F_{\overline{\omega}^*})\geq(1+\underline{b}^{-2})^{-1}\underline{b}^2\max\{|F_{\overline{\sigma}}|,|F_{\overline{\omega}}|\}
\\&\geq&(1+\underline{b}^{-2})^{-1}\underline{b}^2\max\{|F_{\overline{\sigma}^*}|,|F_{\overline{\omega}^*}|\}.
\end{eqnarray*}
If $(\overline{\sigma})_y,(\overline{\omega})_y$ are comparable, then by Lemmas \ref{zz0}, \ref{zsg1}, $\sigma_L,\omega_L$ are incomparable. Thus, from Lemma \ref{ss4}, we have
\[
d(F_{\overline{\sigma}^*},F_{\overline{\omega}^*})\geq d(F_{\overline{\sigma}},F_{\overline{\omega}})\geq2^{-\frac{1}{2}}\underline{a}^2\max\{|F_{\overline{\sigma}^*}|,|F_{\overline{\omega}^*}|\}.
\]
Note that $\underline{a}\leq\underline{b}$. The proof of the lemma is complete.
\end{proof}
\begin{remark}\label{zz4}{\rm
(1) From the definitions of $\mathcal{B}_{n,r}$ and $\mathcal{F}_{n,r}$, one can easily see that, ${\rm card}(\mathcal{F}_{n,r})={\rm card}(\mathcal{B}_{n,r})=\varphi_{n,r}$. (2) Let $K_{n,r}:=\bigcup_{\overline{\sigma}^*\in\mathcal{F}_{n,r}}F_{\overline{\sigma}^*}$. By the definition of $\overline{\sigma}^*$, $|\overline{\sigma}^*_L|-|\sigma_L|=2+2A_2$. Thus, by Lemmas \ref{s4} and \ref{zz2},
\begin{eqnarray*}
\mu(K_{n,r})=\sum_{\sigma\in\Lambda_{n,r}}\mu(F_{\overline{\sigma}^*})\geq \underline{p}^{2(A_2+1)(1+A_1)}\sum_{\sigma\in\Lambda_{n,r}}\mu(F_{\sigma})\geq\underline{p}^{8A_1A_2}.
\end{eqnarray*}
(3) As we showed in \cite[Lemma 4]{KZ:15}, there exists a positive number $D_L$, such that, for every $\alpha_L\subset\mathbb{R}^2$ with cardinality $L$, the following holds:
\[
\int_{F_{\overline{\sigma}^*}}d(x,\alpha_L)^rd\mu(x)\geq D_L\mathcal{E}_r(\overline{\sigma}^*).
\]}
\end{remark}
With the above preparations, we can now establish a characterization for the quantization error by applying \cite[Lemma 3]{KZ:15}.

\begin{proposition}\label{characterization}
We have $e_{\varphi_{n,r},r}^r(\mu)\asymp\sum_{\sigma\in\Lambda_{n,r}}\mathcal{E}_r(\sigma)$.
\end{proposition}
\begin{proof}
For every $\sigma\in\Lambda_{n,r}$, let $C_\sigma$ be an arbitrary point in $F_\sigma$. We have
\[
e_{\varphi_{n,r},r}^r(\mu)\leq\sum_{\sigma\in\Lambda_{n,r}}\mu(F_\sigma)|F_\sigma|^r\lesssim\sum_{\sigma\in\Lambda_{n,r}}\mathcal{E}_r(\sigma).
\]
Let $\mathcal{F}_{n,r}$ be as defined in (\ref{Fjr}). For distinct words $\overline{\sigma}^*,\overline{\omega}^*\in\mathcal{F}_{n,r}$, we have
\begin{equation}\label{gg2}
\mathcal{E}_r(\overline{\sigma}^*)\geq\underline{\eta}_r^{2(A_2+1)}\mathcal{E}_r(\sigma)
\geq\underline{\eta}_r^{2(A_2+1)+1}\mathcal{E}_r(\omega)\geq\underline{\eta}_r^{2A_2+3}\mathcal{E}_r(\overline{\omega}^*).
\end{equation}
 By (\ref{gg2}), Lemma \ref{zz2} and Remark \ref{zz4} (3), the assumptions in Lemma 3 of \cite{KZ:15} are fulfilled for the measure $\mu_{n,r}:=\mu(\cdot|K_{n,r})$. It follows that
\[
e^r_{\varphi_{n,r},r}(\mu)\geq\mu(K_{n,r})e^r_{\varphi_{n,r},r}(\mu_{n,r})
\gtrsim\sum_{\sigma\in\mathcal{F}_{n,r}}\mathcal{E}_r(\overline{\sigma}^*)\gtrsim\sum_{\sigma\in\Lambda_{n,r}}\mathcal{E}_r(\sigma).
\]
This completes the proof of the proposition.
\end{proof}

\section{Auxiliary coding space and auxiliary measures}

\subsection{Auxiliary coding space} Let $G,G_y$ be endowed with discrete topology and let $G^{\mathbb{N}},G_y^{\mathbb{N}}$
be endowed with product topology. Then $G^{\mathbb{N}},G_y^{\mathbb{N}}$ are both metrizable. The corresponding product metric is compatible with the product topology on $\Phi_\infty=G^{\mathbb{N}}\times G_y^{\mathbb{N}}$.
Thus, $\Phi_\infty$ is a compact metric space.

With the next lemma, we show that, if $\sigma,\omega\in\Phi^*$ and $[\sigma]\cap[\omega]\neq\emptyset$, we have
either $[\sigma]\subset[\omega]$, or $[\omega]\subset[\sigma]$. The proof of this lemma is different from that for Lemma \ref{tem1}, because $\Psi^*,\Phi^*$ are endowed with different partial orders.
\begin{lemma}\label{t1}
(1) Let $\sigma,\omega\in\Phi^*$. Assume that $\sigma_L\preceq\omega_L$ and $\sigma_R,\omega_R$ are comparable, then we have $\sigma_R\preceq\omega_R$. (2) For every pair $\sigma,\omega\in\Phi^*$, we have either $[\sigma]\cap[\omega]=\emptyset$, or $[\sigma]\subset[\omega]$, or $[\omega]\subset[\sigma]$.
\end{lemma}
\begin{proof}
(1) Assume that $\sigma_L\preceq\omega_L, \omega_R\precneqq\sigma_R$. We write
\begin{eqnarray*}
&\sigma=((i_1,j_1),\ldots,(i_l,j_l))\times(j_{l+1},\ldots, j_k);\\
&\omega=((i_1,j_1),\ldots,(i_l,j_l),\ldots, (\hat{i}_{l+p_1},\hat{j}_{l+p_1}))\times(j_{l+1},\ldots, j_{k-p_2}),
\end{eqnarray*}
where $p_1\geq 0$ and $p_2\geq 1$. We have, $b_{\sigma_y^-}\geq a_{\sigma_L}>b_{\sigma_y}$. If $p_1=0$, we have
\begin{eqnarray*}
b_{\omega_y}&=&\prod_{h=1}^{l}b_{j_h}\prod_{h=l+1}^{k-p_2}b_{j_h}=\prod_{h=1}^{k-p_2}b_{j_h}\geq \prod_{h=1}^la_{i_hj_h}=a_{\omega_L}.
\end{eqnarray*}
This contradicts the fact that $\omega\in\Phi_l$. Next, we assume that $p_1\geq 1$. Note that $a_{ij}<b_j$ for every $(i,j)\in G$. It follows that
\begin{eqnarray*}
b_{\omega_y}&=&\prod_{h=1}^{l}b_{j_h}\prod_{h=l+1}^{l+p_1}b_{\hat{j}_h}\prod_{h=l+1}^{k-p_2}b_{j_h}\\
&=&\prod_{h=1}^{k-p_2}b_{j_h}\prod_{h=l+1}^{l+p_1}b_{\hat{j}_h}\geq \prod_{h=1}^la_{i_hj_h}\prod_{h=l+1}^{l+p_1}b_{\hat{j}_h}\\
&>&\prod_{h=1}^la_{i_hj_h}\prod_{h=l+1}^{l+p_1}a_{\hat{j}_h\hat{j}_h}.
\end{eqnarray*}
This contradicts the fact that $\omega\in\Phi_{l+p_1}$. Hence, $|\omega_R|\geq|\sigma_R|$ and $\sigma_R\preceq\omega_R$.

(2) If both $\sigma_L,\omega_L$, and $\sigma_R,\omega_R$, are comparable, then by (1), we have $[\sigma]\subset[\omega]$, or $[\omega]\subset[\sigma]$. Otherwise, either $\sigma_L,\omega_L$, or $\sigma_R,\omega_R$, are incomparable, and then we have $[\sigma]\cap[\omega]=\emptyset$.
\end{proof}
\begin{remark}\label{r2}{\rm Based on Lemma \ref{t1}, we obtain the following useful facts.
\begin{enumerate}
\item[(r0)]Let $\sigma=((i_1,j_1),\ldots,(i_l,j_l))\times(j_{l+1},\ldots, j_k)\in\Phi_l$ and $\omega\in\Phi_{l+1}$ with $\sigma\preceq\omega$. Then for some $p\geq 0$ and $(i,j)\in G$, we have
    \[
    \omega=((i_1,j_1),\ldots,(i_l,j_l),(i,j))\times(j_{l+1},\ldots, j_k,\ldots,j_{k+p}).
    \]
\item[(r1)] For every pair of distinct words $\sigma,\omega\in\Phi_l$, we have that $[\sigma]\cap[\omega]=\emptyset$. In fact, if $\sigma_L\neq\omega_L$, then we certainly have that $[\sigma]\cap[\omega]=\emptyset$, since $|\sigma_L|=|\omega_L|=l$; if $\sigma_L=\omega_L$, then $\sigma_R,\omega_R$ are incomparable, and we again have $[\sigma]\cap[\omega]=\emptyset$.

\item[(r2)] For every $\omega\in G^l$, we define $\Omega(\omega):=\{\tau\in G_y^*:\omega\times\tau\in\Phi_l\}$. We have
\[
\Omega(\omega):=\{\tau\in G_y^*:b_{\tau^-}\geq\frac{a_{\omega}}{b_{\omega_y}}>b_\tau\}.
\]
Hence, $G_y^{\mathbb{N}}$ is the disjoint union of the sets $[\tau],\tau\in\Omega(\omega)$. Therefore,
\[
\bigcup_{\sigma\in\Phi_l}[\sigma]=\bigcup_{\omega\in G^l}\bigcup_{\tau\in\Omega(\omega)}[\omega\times\tau]=\Phi_\infty.
\]
\end{enumerate}}
\end{remark}
\begin{remark}\label{r3} {\rm The following facts will also be useful (cf. \cite{Zhu:18}).
\begin{enumerate}
\item[(r3)] It can happen that for some $\sigma=\sigma_L\ast\sigma_R,\omega=\omega_L\ast\omega_R\in\Psi^*$ with $F_{\sigma}^\circ\cap F_\omega^\circ=\emptyset$,
but $[\mathcal{L}(\sigma)]\supset[\mathcal{L}(\omega)]$. This can be seen by considering
  \begin{eqnarray*}
  &&\sigma_L=((i_1,j_1),\ldots,(i_l,j_l)),\;\sigma_R=(j_{l+1},\ldots, j_k);\;j_{l+1}\neq \hat{j}_{l+1};\\
 &&\omega_L=((i_1,j_1),\ldots,(i_l,j_l),(\hat{i}_{l+1},\hat{j}_{l+1})),\;\omega_R=(j_{l+1},\ldots, j_k,\ldots,j_{k+p}).
  \end{eqnarray*}

\item[(r4)]  It can happen that for some $\sigma=\sigma_L\times\sigma_R,\omega=\omega_L\times\omega_R\in\Phi^*$, $[\sigma]\cap[\omega]=\emptyset$, but $F_{\mathcal{L}^{-1}(\omega)}\subset F_{\mathcal{L}^{-1}(\sigma)}$. This can be seen by considering
  \begin{eqnarray*}
  &&\sigma_L=((i_1,j_1),\ldots,(i_l,j_l)),\;\sigma_R=(j_{l+1},\ldots, j_k);\;j_{l+1}\neq j_{l+2};\\
 &&\omega_L=((i_1,j_1),\ldots,(i_l,j_l),(i_{l+1},j_{l+1})),\;\omega_R=(j_{l+2},\ldots, j_k,\ldots,j_{k+p}).
  \end{eqnarray*}
  \end{enumerate}}
\end{remark}

\subsection{Auxiliary measures}
For $\sigma\in\Psi^*$, let $\mathcal{E}_r(\sigma)$ be as defined in (\ref{energy}). We define
\begin{eqnarray*}
&&\mathcal{E}_r(\sigma):=\mathcal{E}_r(\mathcal{L}^{-1}(\sigma)),\;\sigma\in\Phi^*;\\
&&\Lambda_h(\sigma):=\{\rho\in\Phi_{l+h}:\sigma\preceq\rho\};\;I_{h,r}(t):=\sum_{\omega\in\Phi_h}\mathcal{E}_r(\omega)^t
,\;t\geq0,\;h\geq 1.
\end{eqnarray*}
In order to construct an auxiliary measure on $\Phi_\infty$, we need to prove
\begin{equation}\label{sgg2}
\sum_{\rho\in\Lambda_h(\sigma)}\mathcal{E}_r(\rho)^t\asymp\mathcal{E}_r(\sigma)^tI_{h,r}(t).
\end{equation}
For $t=0$, (\ref{sgg2}) trivially becomes ${\rm card}(\Lambda_h(\sigma))\asymp{\rm card}(\Phi_h)$. In the following, we divide the proof of (\ref{sgg2}) into three lemmas.
\begin{lemma}\label{zzz2}
Let $A_6:=[A_4^{-1}]$ and $\sigma\in\Phi^*$. (1) for every $h>A_6$ and $\rho\in\Lambda_h(\sigma)$, we have $|\rho_R|\geq |\sigma_R|+1$; (2) for every $\rho\in\Lambda_1(\sigma)$, we have $|\rho_R|-|\sigma_R|\leq A_1$.
\end{lemma}
\begin{proof}
Let $\sigma=((i_1,j_1),\ldots,(i_l,j_l))\times(j_{l+1},\ldots,j_k)\in\Phi_l$. Assume that $h>A_6$ and $\rho\in\Lambda_h(\sigma)$. We write $\rho_L:=\sigma_L\ast((\hat{i}_1,\hat{j}_1),\ldots,(\hat{i}_h,\hat{j}_h))$. Then we have
\begin{eqnarray*}
a_{\rho_L}&=&\prod_{p=1}^la_{i_pj_p}\prod_{p=1}^ha_{\hat{i}_p\hat{j}_p}\leq \prod_{p=1}^{k-1}b_{j_p}\prod_{p=1}^ha_{\hat{i}_p\hat{j}_p}\\
&\leq&\underline{b}^{-1}\prod_{p=1}^kb_{j_p}\prod_{p=1}^ha_{\hat{i}_p\hat{j}_p}=\underline{b}^{-1}\prod_{p=1}^kb_{j_p}\prod_{p=1}^hb_{\hat{j}_p}
\prod_{p=1}^h\frac{a_{\hat{i}_p\hat{j}_p}}{b_{\hat{j}_p}}\\
&\leq&\underline{b}^{-1}A_3^h\prod_{p=1}^kb_{j_p}\prod_{p=1}^hb_{\hat{j}_p}\;\;({\rm by}\; (\ref{temp2}))\\
&<&\prod_{p=1}^kb_{j_p}\prod_{p=1}^hb_{\hat{j}_p}
\end{eqnarray*}
It follows that $|\rho_R|>k-l$ and (1) follows. (2) can be proved similarly.
\end{proof}
\begin{lemma}\label{lem1b}
For every $t\geq 0$, there exists a number $h_{1,r}(t)>0$ such that for every $\sigma\in\Phi^*$ and $h\geq 1$, the following holds:
\begin{equation}\label{sss2}
\sum_{\rho\in\Lambda_h(\sigma)}\mathcal{E}_r(\rho)^t\leq h_{1,r}(t)\mathcal{E}_r(\sigma)^tI_{h,r}(t).
\end{equation}
\end{lemma}
\begin{proof}
Let $N:={\rm card}(G)$. For $\sigma\in\Phi^*$ and $h\geq 1$, by Lemma \ref{zzz2},, we have ${\rm card}(\Lambda_h(\sigma))\leq N^hm^{hA_1}$. Hence,
\[
\sum_{\rho\in\Lambda_h(\sigma)}\mathcal{E}_r(\rho)^t\leq N^hm^{hA_1}(\overline{p}\overline{a}^r)^h\mathcal{E}_r(\sigma)^t< N^hm^{hA_1}\overline{\eta}_r^{ht}\mathcal{E}_r(\sigma)^t.
\]
Let $\xi_{1,r}(t):=\max\limits_{1\leq h\leq A_6}(I_{h,r}(t)^{-1}N^hm^{hA_1}\overline{\eta}_r^{ht})$. For every $1\leq h\leq A_6$, we have
\[
\sum_{\rho\in\Lambda_h(\sigma)}\mathcal{E}_r(\rho)^t\leq \xi_{1,r}(t)\mathcal{E}_r(\sigma)^tI_{h,r}(t).
\]
Next, we assume that $h>A_6$. Let $\sigma=((i_1,j_1),\ldots,(i_l,j_l))\times(j_{l+1},\ldots,j_k)$ be given.
Let $\rho$ be an arbitrary word in $\Lambda_h(\sigma)$. By Lemma \ref{zzz2}, we have $|\rho_R|\geq|\sigma_R|+1$. Write
\begin{equation*}
\rho=((i_1,j_1),\ldots,(i_l,j_l),(\hat{i}_1,\hat{j}_1),\ldots,(\hat{i}_h,\hat{j}_h))\times(j_{l+1},\ldots,j_k,\hat{j}_{h+1},\ldots, \hat{j}_{\hat{k}}).
\end{equation*}
By the definition of $\Phi^*$, we have
\begin{eqnarray}
&\prod_{p=1}^{k-1}b_{j_p}\geq \prod_{p=1}^la_{i_pj_p}>\prod_{p=1}^kb_{j_p};\nonumber
\\&\prod_{p=1}^kb_{j_p}\prod_{p=1}^{\hat{k}-1}b_{\hat{j}_p}\geq \prod_{p=1}^la_{i_pj_p}\prod_{p=1}^ha_{\hat{i}_p\hat{j}_p}>\prod_{p=1}^kb_{j_p}\prod_{p=1}^{\hat{k}}b_{\hat{j}_p}.\label{ggg4}
\end{eqnarray}
As a consequence, we obtain
\[
\prod_{p=1}^{\hat{k}-1}b_{\hat{j}_p}>\prod_{p=1}^ha_{\hat{i}_p\hat{j}_p}>b_{j_k}\prod_{p=1}^{\hat{k}}b_{\hat{j}_p}.
\]
We distinguish between the following two cases.
\begin{enumerate}
\item[(i)] $\prod_{p=1}^ha_{\hat{i}_p\hat{j}_p}>\prod_{p=1}^{\hat{k}}b_{\hat{j}_p}$. In this case, we define
\[
\omega(\rho):=((\hat{i}_1,\hat{j}_1),\ldots,(\hat{i}_h,\hat{j}_h))\times(\hat{j}_{h+1},\ldots, \hat{j}_{\hat{k}}).
\]
Then $\omega(\rho)\in\Phi_h$. Let $\widetilde{\Phi}_{h,1}$ denote the set of all such words $\omega(\rho)$ and let $\widetilde{\Lambda}_{h,1}(\sigma)$ denote the set of the words $\rho$ in this case.
\item[(ii)]$\prod_{p=1}^ha_{\hat{i}_p\hat{j}_p}\leq\prod_{p=1}^{\hat{k}}b_{\hat{j}_p}$. In this case, we define
\[
\omega(\rho):=((\hat{i}_1,\hat{j}_1),\ldots,(\hat{i}_h,\hat{j}_h))\times(\hat{j}_{h+1},\ldots, \hat{j}_{\hat{k}},j_k).
\]
One can see that $\omega(\rho)\in\Phi_h$. Let $\widetilde{\Phi}_{h,2}$ denote the set of all such words $\omega$ and let $\widetilde{\Lambda}_{h,2}(\sigma)$ the set of the words $\rho$ in this case.
\end{enumerate}

We clearly have that, $\widetilde{\Lambda}_{h,1}(\sigma)\cap\widetilde{\Lambda}_{h,2}(\sigma)=\emptyset,\Lambda_h(\sigma)=\widetilde{\Lambda}_{h,1}(\sigma)\cup\widetilde{\Lambda}_{h,2}(\sigma)$. Further, we have that $\widetilde{\Phi}_{h,1}\cap\widetilde{\Phi}_{h,2}=\emptyset$. Otherwise, there would exist some $\rho^{(1)}\in\widetilde{\Lambda}_{h,1}(\sigma)$ and $\rho^{(2)}\in\widetilde{\Lambda}_{h,2}(\sigma)$ such that $\omega(\rho^{(1)})=\omega(\rho^{(2)})$. Then $(\rho^{(1)})_L=(\rho^{(2)})_L$ and $(\rho^{(2)})_R=((\rho^{(1)})_R)^-$, contradicting (\ref{ggg4}). For $t=0$, we simply have that ${\rm card}(\Lambda_h(\sigma))\leq {\rm card}(\Phi_h)$. For $t>0$, we have
\begin{eqnarray*}
I_{h,r}(t)&\geq&\sum_{\omega\in\widetilde{\Phi}_{h,1}}\mathcal{E}_r(\omega)^t+\sum_{\omega\in\widetilde{\Phi}_{h,2}}\mathcal{E}_r(\omega)^t\\
&\geq&\sum_{\rho\in\widetilde{\Lambda}_{h,1}(\sigma)}\frac{\mathcal{E}_r(\rho)^t}{\mathcal{E}_r(\sigma)^t}+\underline{q}^t\sum_{\rho\in\widetilde{\Lambda}_{h,2}(\sigma)}\frac{\mathcal{E}_r(\rho)^t}{\mathcal{E}_r(\sigma)^t}\\
&\geq&\underline{q}^t\mathcal{E}_r(\sigma)^{-t}\sum_{\rho\in\Lambda_h(\sigma)}\mathcal{E}_r(\rho)^t.
\end{eqnarray*}
It is sufficient to define $h_{1,r}(t):=\max\{\xi_{1,r}(t),\underline{q}^{-t}\}$.
\end{proof}

\begin{lemma}\label{lem1a}
For every $t\geq 0$, there exists a number $h_{2,r}(t)>0$, such that for every $\sigma\in\Phi^*$ and $h\geq 1$, the following holds:
\begin{equation}\label{d1}
\sum_{\rho\in\Lambda_h(\sigma)}\mathcal{E}_r(\rho)^t\geq h_{2,r}(t)\mathcal{E}_r(\sigma)^tI_{h,r}(t).
\end{equation}
\end{lemma}
\begin{proof}
For $\sigma\in\Phi^*$ and $h\geq 1$, we have ${\rm card}(\Lambda_h(\sigma))\geq N^h$. By Lemma \ref{zzz2},
\[
\sum_{\rho\in\Lambda_h(\sigma)}\mathcal{E}_r(\rho)^t\geq N^h\underline{\eta}_r^{ht}\mathcal{E}_r(\sigma)^t.
\]
Let $A_7:=[\frac{\log\underline{b}}{\log\overline {b}}]+2$ and  $A_8:=[A_7A_4^{-1}]+1$. We define
 \[
 \xi_{2,r}(t):=\min\limits_{1\leq h\leq A_8}(I_{h,r}(t)^{-1}N^h\underline{\eta}_r^{ht}).
 \]
 Then for every $1\leq h\leq A_8$, we have
\[
\sum_{\rho\in\Lambda_h(\sigma)}\mathcal{E}_r(\rho)^t\geq \xi_{2,r}(t)\mathcal{E}_r(\sigma)^tI_{h,r}(t).
 \]
 Next, we assume that $h>A_8$. Let $\sigma=((i_1,j_1),\ldots,(i_l,j_l))\times(j_{l+1},\ldots,j_k)$ be given. Let $\omega$ be an arbitrary word in $\Phi_h$. By Lemma \ref{zzz1}, we know that $|\omega_R|>A_7$. We write $\omega=((\hat{i}_1,\hat{j}_1),\ldots,(\hat{i}_h,\hat{j}_h))\times (\hat{j}_{h+1},\ldots, \hat{j}_{\hat{k}})$.  We have
\begin{eqnarray}\label{ss1}
\prod_{p=1}^{k-1}b_{j_p}\geq \prod_{p=1}^la_{i_pj_p}>\prod_{p=1}^kb_{j_p};\;\;\prod_{p=1}^{\hat{k}-1}b_{\hat{j}_p}
\geq \prod_{p=1}^ha_{\hat{i}_p\hat{j}_p}>\prod_{p=1}^{\hat{k}}b_{\hat{j}_p}.
\end{eqnarray}
It follows that
\[
\prod_{p=1}^{k-1}b_{j_p}\prod_{p=1}^{\hat{k}-1}b_{\hat{j}_p}\geq \prod_{p=1}^la_{i_pj_p}\prod_{p=1}^ha_{\hat{i}_p\hat{j}_p}>\prod_{p=1}^kb_{j_p}\prod_{p=1}^{\hat{k}}b_{\hat{j}_p}.
\]
We need to distinguish between the following two cases:
\begin{enumerate}
\item[(1)] $\prod_{p=1}^kb_{j_p}\prod_{p=1}^{\hat{k}-1}b_{\hat{j}_p}\geq \prod_{p=1}^la_{i_pj_p}\prod_{p=1}^ha_{\hat{i}_p\hat{j}_p}$. In this case, we define
\[
\rho(\omega):=(\sigma_L\ast\omega_L)\times(\sigma_R\ast\omega_R).
\]
We have that $\rho(\omega)\in\Lambda_h(\sigma)$. We denote by $\widehat{\Lambda}_{h,1}(\sigma)$ the of set such words $\rho(\omega)$ and denote the set of the words $\omega$ by $\widehat{\Phi}_{h,1}$.

\item[(2)] $\prod_{p=1}^kb_{j_p}\prod_{p=1}^{\hat{k}-1}b_{\hat{j}_p}<\prod_{p=1}^la_{i_pj_p}\prod_{p=1}^ha_{\hat{i}_p\hat{j}_p}$. We have
\begin{eqnarray}\label{temp3}
\prod_{p=1}^kb_{j_p}\prod_{p=1}^{\hat{k}-A_7}b_{\hat{j}_p}&=&\prod_{p=1}^{k-1}b_{j_p}\prod_{p=1}^{\hat{k}-1}b_{\hat{j}_p}
\frac{b_{j_k}}{\prod_{p=\hat{k}-A_7+1}^{\hat{k}-1}b_{\hat{j}_p}}\nonumber\\&\geq&\prod_{p=1}^{k-1}b_{j_p}\prod_{p=1}^{\hat{k}-1}b_{\hat{j}_p}\frac{\underline{b}}{\overline{b}^{A_7-1}}
\nonumber\\&>&\prod_{p=1}^{k-1}b_{j_p}\prod_{p=1}^{\hat{k}-1}b_{\hat{j}_p}\geq\prod_{p=1}^la_{i_pj_p}\prod_{p=1}^ha_{\hat{i}_p\hat{j}_p}.
\end{eqnarray}
Hence, there exists a unique integer $2\leq H\leq A_7$, such that
\[
\prod_{p=1}^kb_{j_p}\prod_{p=1}^{\hat{k}-H}b_{\hat{j}_p}\geq\prod_{p=1}^la_{i_pj_p}\prod_{p=1}^ha_{\hat{i}_p\hat{j}_p}>\prod_{p=1}^kb_{j_p}\prod_{p=1}^{\hat{k}-H+1}b_{\hat{j}_p}.
\]
 Let $\omega_R|_{\hat{k}-h-H+1}:=(\hat{j}_{h+1},\ldots,\hat{j}_{\hat{k}-H+1})$. We define
 \[
 \rho(\omega):=(\sigma_L\ast\omega_L)\times(\sigma_R\ast(\omega_R|_{\hat{k}-h-H+1}).
  \]
 Then $\rho(\omega)\in\Lambda_h(\sigma)$. We denote the set of such words $\rho(\omega)$ by $\widehat{\Lambda}_{h,2}(\sigma)$ and denote the set of the words $\omega$ by $\widehat{\Phi}_{h,2}$.
\end{enumerate}

We need to observe the following facts.

(i) We clearly have that, $\widehat{\Phi}_{h,1}\cap\widehat{\Phi}_{h,2}=\emptyset,\Phi_h=\widehat{\Phi}_{h,1}\cup\widehat{\Phi}_{h,2}$. Also, we have, $\widehat{\Lambda}_{h,1}(\sigma)\cap\widehat{\Lambda}_{h,2}(\sigma)=\emptyset$. In fact, a word $\rho(\omega)\in\widehat{\Lambda}_{h,1}(\sigma)$ can not be obtained by any $\tau\in\widehat{\Phi}_{h,2}$
and vice versa. Otherwise, we would have $\omega_L=\tau_L$ and $\tau_R\precneqq\omega_R$, contradicting the fact that $\omega,\tau\in\Phi_h$.

(ii) For different words $\omega,\tau\in \widehat{\Phi}_{h,1}$, we have $\rho(\omega)\neq\rho(\tau)$.

(iii) There exist at most $\sum_{p=1}^{A_7}m^h$ words in $\widehat{\Phi}_{h,2}$ that determine the same word in $\widehat{\Lambda}_{h,2}(\sigma)$, because of the absence of $\hat{j}_{\hat{p}}, \hat{k}-H+2\leq p\leq \hat{k}$ in $\rho(\omega)$. To see this, we fix an $\omega\in\widehat{\Phi}_{h,2}$. For $\tau\in G_y^*$, let $\hat{\omega}^{(\tau)}:=\omega_L\times (\omega_R|_{\hat{k}-h-H+1}\ast \tau)$. Let $\rho^{-1}(\rho(\omega)):=\{\hat{\omega}\in\Phi_{h,2}:\rho(\hat{\omega})=\rho(\omega)\}$. We have
\[
\rho^{-1}(\rho(\omega))=:\langle\omega\rangle=\big\{\hat{\omega}^{(\tau)}:\tau\in G_y^*\big\}\cap\Phi_{h,2}.
 \]
 For every $\hat{\omega}^{(\tau)}\in \langle\omega\rangle$, by (\ref{temp3}), we have $|\tau|<A_7$. It follows that
    \begin{equation}\label{angle}
    \langle\omega\rangle\subset\bigg\{\hat{\omega}^{(\tau)}:\tau\in \bigcup_{h=1}^{A_7}G_y^h\bigg\}.
    \end{equation}

For $\omega\in\widehat{\Phi}_{h,2}$, we take an arbitrary word of $\langle\omega\rangle$ and denote the set of such words by $\widehat{\Phi}^\flat_{h,2}$. By (\ref{angle}),
\begin{equation}\label{angle2}
\sum_{\omega\in\widehat{\Phi}_{h,2}}\mathcal{E}_r(\omega)^t=\sum_{\omega\in\widehat{\Phi}_{h,2}^\flat}\sum_{\widetilde{\omega}\in\langle\omega\rangle}\mathcal{E}_r(\widetilde{\omega})^t
\leq\sum_{\omega\in\widehat{\Phi}_{h,2}^\flat}\frac{\mathcal{E}_r(\rho(\omega))^t}{\mathcal{E}_r(\sigma)^t}\sum_{p=1}^{A_7}\big(\sum_{j\in G_y}q_j^t\big)^p.
\end{equation}
Let $\xi_3(t):=\sum_{p=1}^{A_7}\big(\sum_{j\in G_y}q_j^t\big)^p$. Using (\ref{angle2}), we deduce
\begin{eqnarray*}
\sum_{\rho\in\Lambda_h(\sigma)}\mathcal{E}_r(\rho)^t&\geq&\sum_{\rho\in\widehat{\Lambda}_{h,1}(\sigma)}\mathcal{E}_r(\rho)^t
+\sum_{\rho\in\widehat{\Lambda}_{h,2}(\sigma)}\mathcal{E}_r(\rho)^t\\
&=&\mathcal{E}_r(\sigma)^t\sum_{\omega\in\widehat{\Phi}_{h,1}}\mathcal{E}_r(\omega)^t
+\sum_{\omega\in\widehat{\Phi}^\flat_{h,2}}\mathcal{E}_r(\rho(\omega))^t\\&\geq&
\mathcal{E}_r(\sigma)^t\sum_{\omega\in\widehat{\Phi}_{h,1}}\mathcal{E}_r(\omega)^t
+\mathcal{E}_r(\sigma)^t\xi_3(t)^{-1}\sum_{\omega\in\widehat{\Phi}_{h,2}}\mathcal{E}_r(\omega)^t\\
&\geq&\min\bigg\{1,\xi_3(t)^{-1}\bigg\}\mathcal{E}_r(\sigma)^t I_{h,r}(t).
\end{eqnarray*}
Thus, (\ref{d1}) is fulfilled by defining $h_{2,r}(t):=\min\{1,\xi_3(t)^{-1},\xi_{2,r}(t)\}$.
\end{proof}

 Combining Lemmas \ref{lem1b}, \ref{lem1a}, we obtain
\begin{lemma}\label{z4}
For every $t\geq 0$ and $k,p\in\mathbb{N}$, we have
\[
h_{2,r}(t) I_{k,r}(t)I_{p,r}(t)\leq I_{k+p,r}(t)\leq h_{1,r}(t) I_{k,r}(t)I_{p,r}(t).
\]
\end{lemma}
\begin{proof}
Using Lemmas \ref{lem1b}, \ref{lem1a}, we have
\begin{eqnarray*}
I_{k+p,r}(t)=\sum_{\sigma\in\Phi_k}\sum_{\rho\in\Lambda_p(\sigma)}\mathcal{E}_r(\rho)^t\geq h_{2,r}(t)I_{k,r}(t)I_{p,r}(t).
\end{eqnarray*}
The second inequality can be obtained similarly.
\end{proof}

Using Lemma \ref{z4}, we can obtain  the following standard result.
\begin{proposition}\label{sss3}
(1) For every $t\geq 0$, $\lim_{k\to\infty}\frac{1}{k}\log I_{k,r}(t)=:g_r(t)$ exists. (2) There exists a unique $s_r>0$ such that for $t_r=\frac{s_r}{s_r+r}$, we have $g_r(t_r)=0$.
(3) For $C(t):=h_{2,r}(t)^{-1}h_{1,r}(t)$ and $k,p\geq 1$, we have,
\[
C(t)^{-1}I_{k,r}(t_r)\leq I_{p,r}(t_r)\leq C(t) I_{k,r}(t_r).
\]
\end{proposition}
\begin{proof}
(1) is an easy consequence of Lemma \ref{z4} and \cite[Corollary 1.2]{Fal:97}. (2) Along the line of \cite[Lemma 5.2]{Fal:97}, one can see that $g_r$ is strictly decreasing and  continuous. Note that $g_r(0)\geq\log N>0, g_r(1)\leq r\log\overline{a}<0$. Thus, (2) follows from the continuity of $g_r$. (3) For every $k\geq 1$, we have (cf. \cite[(5.10)]{Fal:97})
\begin{eqnarray*}
\frac{1}{k}(\log I_{k,r}(t_r)+\log h_{1,r}(t_r))\geq\inf_{p\geq 1}\frac{1}{p}(\log I_{p,r}(t_r)+\log h_{1,r}(t_r))=0.
\end{eqnarray*}
It follows that $I_{k,r}(t_r)\geq h_{1,r}(t_r)^{-1}$. Similarly, we have $I_{k,r}(t_r)\leq h_{2,r}(t_r)^{-1}$. This completes the proof of (3).
\end{proof}

With the help of Lemmas \ref{lem1b}, \ref{lem1a} and Proposition \ref{sss3}, we are now able to construct an auxiliary probability measure on $\Phi_\infty$ by applying Prohorov's theorem. Recall that a family $\mathcal{\pi}$ of probability measures on a metric space $X$, is said to be \emph{tight} if for every $\epsilon\in(0,1)$, there exists some compact subset $K$ of $X$ such that $\inf_{\nu\in\mathcal{\pi}}\nu(K)\geq 1-\epsilon$. In particular, if $X$ is a compact metric space, then every family $\mathcal{\pi}$ of probability measures on $X$ is tight.
\begin{lemma}\label{z2}
There exists a Borel probability measure $\lambda$ on $\Phi_\infty$ such that
\[
\lambda([\sigma])\asymp\mathcal{E}_r(\sigma)^{t_r},\;\sigma\in\Phi^*.
\]
\end{lemma}
\begin{proof}
For every $k\geq 1$ and $\sigma\in\Phi_k$, let $C_\sigma$ be an arbitrary point in $[\sigma]$ and let $\delta_{C_\sigma}$ denote the Dirac measure at the point $C_\sigma$. We define
\[
\lambda_k=\frac{1}{I_{k,r}(t_r)}\sum_{\sigma\in\Phi_k}\mathcal{E}_r(\sigma)^{t_r}\delta_{C_\sigma}.
\]
Note that $\Phi_\infty$ is a compact metric space; so $(\lambda_k)_{k=1}^\infty$ is tight.
By Prohorov's Theorem (cf. \cite[Theorem 5.1]{Billings:84}), there exist a subsequence $(\lambda_{k_i})_{i=1}^\infty$ and a probability measure $\lambda$ on $\Phi_\infty$
such that $\lambda_{k_i}$ converges weakly to $\lambda$. Let $n\geq 1$ and $\sigma\in\Phi_n$ be given. For every $i>n$, using Lemmas \ref{lem1b}, \ref{lem1a} and Proposition \ref{sss3}, we deduce
\begin{eqnarray}
\lambda_{k_i}([\sigma])=\sum_{\rho\in\Lambda_{k_i-n}(\sigma)}\lambda_{k_i}([\rho])=\frac{1}{I_{k_i,r}(t_r)}\sum_{\rho\in\Lambda_{k_i-n}(\sigma)}\mathcal{E}_r(\rho)^{t_r}
\asymp \mathcal{E}_r(\sigma)^{t_r}.\label{p1}
\end{eqnarray}
Note that $[\sigma_L]$ ($[\sigma_R]$) is both open and closed in $G^\mathbb{N}$ ($G_y^{\mathbb{N}}$). Thus $[\sigma]$ is clopen in $\Phi_\infty$. Because $\Phi_\infty$ is compact,  $[\sigma]$ is also compact. By (\ref{p1}), it follows  that $\lambda([\sigma])\asymp\mathcal{E}_r(\sigma)^{t_r}$.
\end{proof}

\section{Proofs of Theorem \ref{mthm} and Proposition \ref{mthm2}}

\subsection{Proof of Theorems \ref{mthm}}
By \cite[Lemma 3.4]{Zhu:23}, for the proof of Theorem \ref{mthm}, it is sufficient to show
\begin{equation}\label{q1}
0<\liminf_{n\to\infty}\sum_{\sigma\in\Lambda_{n,r}}\mathcal{E}_r(\sigma)^{t_r}\leq\limsup_{n\to\infty}\sum_{\sigma\in\Lambda_{n,r}}\mathcal{E}_r(\sigma)^{t_r}<\infty.
\end{equation}

Next, we give the proof for the first inequality of (\ref{q1}), in an analogous manner to that for \cite[Proposition 3.2]{Zhu:18}.
\begin{lemma}\label{k1}
We have $\sum_{\sigma\in\Lambda_{n,r}}\mathcal{E}_r(\sigma)^{t_r}\lesssim 1$.
\end{lemma}
\begin{proof}
For every $n\geq 1$, let $\Lambda_{n,r}$ be as defined in (\ref{lambdajr}). We define
\[
\hat{\Lambda}_{n,r}:=\{\hat{\sigma}=\mathcal{L}(\sigma):\sigma\in\Lambda_{n,r}\}.
\]
By Lemma \ref{t1}, for $\hat{\sigma},\hat{\omega}\in\hat{\Lambda}_{n,r}$, either the sets $[\hat{\sigma}], [\hat{\omega}]$, are disjoint, or one is contained in the other. Also, by Remark \ref{cardlambdajr} (1), we have
\[
\underline{\eta}_r\mathcal{E}_r(\hat{\omega})\leq\mathcal{E}_r(\hat{\sigma})\leq\underline{\eta}_r^{-1}\mathcal{E}_r(\hat{\omega}),\;\hat{\sigma},\hat{\omega}\in\hat{\Lambda}_{n,r}.
\]
If $[\hat{\omega}]\subset[\hat{\sigma}]$, then as we did in \cite[Lemma 3.1]{Zhu:18}, one can see that, there exists a constant $H_{1,r}\geq 1$ such
that $\big||\hat{\sigma}_L|-|\hat{\omega}_L|\big|\leq H_{1,r}$. This allows us to select a subset $\hat{\Lambda}_{n,r}^\flat$ of $\hat{\Lambda}_{n,r}$ such that, the
sets $[\hat{\sigma}],\hat{\sigma}\in\hat{\Lambda}_{n,r}^\flat$, are pairwise disjoint and
\[
\hat{\Lambda}_{n,r}=\bigcup_{\hat{\sigma}\in\hat{\Lambda}_{n,r}^\flat}\Gamma(\hat{\sigma}),\;{\rm with}\;\;
\Gamma(\hat{\sigma})=\{\hat{\omega}:\hat{\sigma}\preceq\hat{\omega},\hat{\omega}\in\hat{\Lambda}_{n,r}\}.
\]
Combining the preceding equality with Lemma \ref{z2}, we obtain
\begin{eqnarray*}
\sum_{\hat{\sigma}\in\hat{\Lambda}_{n,r}}\mathcal{E}_r(\hat{\sigma})^{t_r}
\asymp\sum_{\sigma\in\hat{\Lambda}_{n,r}}\lambda([\hat{\sigma}])
\leq (H_{1,r}+1)\sum_{\sigma\in\hat{\Lambda}_{n,r}^\flat}\lambda([\hat{\sigma}])\leq H_{1,r}+1.
\end{eqnarray*}
This completes the proof for the lemma.
\end{proof}

In the following, we are going to prove the last inequality in (\ref{q1}). For every $\sigma\in\Phi_1$, we define $\sigma^-:=\theta$. Let $l\geq 2$ and $\sigma=((i_1,j_1),\ldots,(i_l,j_l))\times(j_{l+1},\ldots,j_k)\in\Phi_l$, we have,
\[
\prod_{h=1}^{l-1}a_{i_hj_h}>a_{i_lj_l}^{-1}\prod_{h=1}^kb_{j_h}=\frac{b_{j_l}}{a_{i_lj_l}}\prod_{h=1}^{l-1}b_{j_h}\prod_{h=l+1}^kb_{j_h}
>\prod_{h=1}^{l-1}b_{j_h}\prod_{h=l+1}^kb_{j_h}.
\]
There exists a unique integer $p\geq 0$, such that the following inequalities hold:
\[
\prod_{h=1}^{l-1}b_{j_h}\prod_{h=l+1}^{k-p-1}b_{j_h}\geq\prod_{h=1}^{l-1}a_{i_hj_h}>\prod_{h=1}^{l-1}b_{j_h}\prod_{h=l+1}^{k-p}b_{j_h}.
\]
We define $\sigma^-:=\sigma_L^-\times(j_{l+1},\ldots,j_{k-p})$. By Lemma \ref{zzz2}, we obtain
\begin{equation}\label{sgg1}
\underline{\eta}_r\mathcal{E}_r(\sigma^-)\leq\mathcal{E}_r(\sigma)\leq\overline{p}\;\overline{a}^r\mathcal{E}_r(\sigma^-)<\overline{\eta}_r\mathcal{E}_r(\sigma^-).
\end{equation}
Let $\mathcal{S}_{n,r}:=\{\sigma\in\Psi^*:\underline{\eta}_r^{n+1}\leq\mathcal{E}_r(\sigma)<\underline{\eta}_r^n\}$. We define
\begin{equation}\label{G1}
G_1(\sigma):=\{\omega\in\mathcal{S}_{n,r}:\sigma\preceq\omega\},\;\sigma\in\mathcal{S}_{n,r}.
\end{equation}
\begin{remark}\label{r1}
{\rm Let $ M_r:=[\frac{\log\underline{\eta}_r}{\log\overline{\eta}_r}]$. For every $\omega\in G_1(\sigma)$, we have, $|\omega_L|-|\sigma_L|\leq M_r$. In fact, assume
that $|\omega_L|-|\sigma_L|>M_r$. Then we have
\[
\mathcal{E}_r(\omega)\leq\mathcal{E}_r(\sigma)\cdot\overline{\eta}_r^{M_r}<\underline{\eta}_r^{n+1}.
\]
This contradicts the definition of $G_1(\sigma)$. Thus, for $\sigma\in\mathcal{S}_{n,r}\cap\Psi_l$, we have that $G_1(\sigma)\subset\bigcup_{h=0}^{M_r}\Gamma_h(\sigma)$, where
$\Gamma_h(\sigma):=\{\rho\in\Psi_{l+h}:\sigma\preceq\rho\}$.
}\end{remark}
The following lemma is an analogue of \cite[Lemma 4.1]{Zhu:25}.
\begin{lemma}\label{k3}
There exists a constant $H_{2,r}>0$ such that, for every $\sigma\in\mathcal{S}_{n,r}$,
\[
\sum_{\omega\in G_1(\sigma)}\mathcal{E}_r(\omega)^{t_r}\leq H_{2,r}\mathcal{E}_r(\sigma)^{t_r}.
\]
\end{lemma}
\begin{proof}
By Lemma \ref{s4} (i), for every $h\geq 1$, we have that ${\rm card}(\Gamma_h(\sigma))\leq N^hm^{hA_1}$.
This, along with (\ref{z3}), yields
\begin{eqnarray*}
\sum_{\omega\in G_1(\sigma)}\mathcal{E}_r(\omega)^{t_r}\leq\sum_{h=0}^{M_r}\sum_{\omega\in\Gamma_h(\sigma)}\mathcal{E}_r(\omega)^{t_r}\leq\sum_{h=0}^{M_r} N^hm^{hA_1}\overline{\eta}_r^{ht_r}\mathcal{E}_r(\sigma)^{t_r}.
\end{eqnarray*}
It is sufficient to define $ H_{2,r}:=\sum_{h=0}^{M_r}N^hm^{hA_1}\overline{\eta}_r^{ht_r}$.
\end{proof}

\begin{lemma}\label{k2}
We have $\sum_{\sigma\in\Lambda_{n,r}}\mathcal{E}_r(\sigma)^{t_r}\gtrsim 1$.
\end{lemma}
\begin{proof}
For every $n\geq 1$, we define
\begin{eqnarray*}
\Gamma_{n,r}:=\{\sigma\in\Phi^*:\mathcal{E}_r(\sigma^-)\geq\underline{\eta}_r^n>\mathcal{E}_r(\sigma)\}.
\end{eqnarray*}
By Lemma \ref{t1} and Remark \ref{r2} (r2), for every pair $\sigma,\omega\in\Gamma_{n,r}$, we have
\[
[\sigma]\cap[\omega]=\emptyset,\;\bigcup_{\sigma\in\Gamma_{n,r}}[\sigma]=\Phi_\infty.
\]
From this and Lemma \ref{z2}, it follows that
\begin{eqnarray}\label{t2}
\sum_{\sigma\in\Gamma_{n,r}}\mathcal{E}_r(\sigma)^{t_r}\asymp\sum_{\sigma\in\Gamma_{n,r}}\lambda([\sigma])=1.
\end{eqnarray}
Now we connect the words in $\Gamma_{n,r}$ with the approximate squares. We define
\[
\widetilde{\Gamma}_{n,r}:=\{\widetilde{\sigma}:=\mathcal{L}^{-1}(\sigma)=\sigma_L\ast\sigma_R:\sigma\in\Gamma_{n,r}\}.
\]
From (\ref{sgg1}) and the definition of $\Gamma_{n,r}$, we know that, $\widetilde{\Gamma}_{n,r}\subset\mathcal{S}_{n,r}$; and every $\widetilde{\sigma}\in\widetilde{\Gamma}_{n,r}$ corresponds to an approximate square. We define
\[
G(\widetilde{\sigma})=\{\widetilde{\omega}\in\widetilde{\Gamma}_{n,r}:F_{\widetilde{\omega}}\subset F_{\widetilde{\sigma}}\},\;\;\widetilde{\sigma}\in\widetilde{\Gamma}_{n,r}.
\]
Then $G(\widetilde{\sigma})\subset G_1(\widetilde{\sigma})$. There exists a subset $\widetilde{\Gamma}_{n,r}^\flat$ of $\widetilde{\Gamma}_{n,r}$ such that
\[
F_{\widetilde{\sigma}}^\circ\cap F_{\widetilde{\omega}}^\circ=\emptyset,\;\widetilde{\sigma},\widetilde{\omega}\in\widetilde{\Gamma}_{n,r}^\flat;
\;\; \widetilde{\Gamma}_{n,r}=\bigcup_{\widetilde{\sigma}\in\widetilde{\Gamma}_{n,r}^\flat}G(\widetilde{\sigma}).
\]
Using this, Lemma \ref{k3} and (\ref{t2}), we deduce
\begin{eqnarray}\label{s1}
\sum_{\widetilde{\sigma}\in\widetilde{\Gamma}_{n,r}^\flat}\mathcal{E}_r(\widetilde{\sigma})^{t_r}
\geq H_{2,r}^{-1}\sum_{\widetilde{\sigma}\in\widetilde{\Gamma}_{n,r}}\mathcal{E}_r(\widetilde{\sigma})^{t_r}
=H_{2,r}^{-1}\sum_{\sigma\in\Gamma_{n,r}}\mathcal{E}_r(\sigma)^{t_r}\asymp 1.
\end{eqnarray}

Next, we compare the words in $\widetilde{\Gamma}_{n,r}^\flat$ with those in $\Lambda_{n,r}$. We define
\begin{equation*}
\Lambda_{n,r}^\flat:=\bigg\{\sigma\in\Lambda_{n,r}: F^\circ_\sigma\cap\bigg(\bigcup_{\widetilde{\sigma}\in\widetilde{\Gamma}_{n,r}^\flat}F_{\widetilde{\sigma}}^\circ\bigg)\neq\emptyset\bigg\}.
\end{equation*}
We need to divide $\Lambda_{n,r}^\flat$ and $\widetilde{\Gamma}_{n,r}^\flat$ into two subsets:
\begin{eqnarray*}
&&\Lambda_{n,r}^\flat(1):=\{\sigma\in\Lambda_{n,r}: F_\sigma\subseteq F_{\widetilde{\sigma}}\;{\rm for\;some}\;\widetilde{\sigma}\in\widetilde{\Gamma}_{n,r}^\flat\},\\
&&\Lambda_{n,r}^\flat(2):=\{\sigma\in\Lambda_{n,r}: F_\sigma\supsetneq F_{\widetilde{\sigma}}\;{\rm for\;some}\;\widetilde{\sigma}\in\widetilde{\Gamma}_{n,r}^\flat\},\\
&&\widetilde{\Gamma}_{n,r}^\flat(1):=\{\widetilde{\sigma}\in\widetilde{\Gamma}_{n,r}^\flat: F_\sigma\subseteq F_{\widetilde{\sigma}}\;{\rm for\;some}\;\sigma\in\Lambda_{n,r}^\flat\},\\
&&\widetilde{\Gamma}_{n,r}^\flat(2):=\{\widetilde{\sigma}\in\widetilde{\Gamma}_{n,r}^\flat: F_\sigma\supsetneq F_{\widetilde{\sigma}}\;{\rm for\;some}\;\sigma\in\Lambda_{n,r}^\flat\}.
\end{eqnarray*}
By Lemma \ref{t1}, $\Lambda_{n,r}^\flat$ is the disjoint union of $\Lambda_{n,r}^\flat(1)$ and $\Lambda_{n,r}^\flat(2)$. We define
\begin{eqnarray*}
&&S(\sigma):=\{\widetilde{\sigma}\in\widetilde{\Gamma}_{n,r}^\flat:F_\sigma\supsetneq F_{\widetilde{\sigma}}\};\;\sigma\in\Lambda_{n,r}^\flat(2);\\
&&T(\widetilde{\sigma}):=\{\sigma\in\Lambda_{n,r}:F_\sigma\subseteq F_{\widetilde{\sigma}}\},\;\widetilde{\sigma}\in\widetilde{\Gamma}_{n,r}^\flat(1).
\end{eqnarray*}
We have $\Lambda_{n,r}^\flat(1)=\bigcup_{\widetilde{\sigma}\in\widetilde{\Gamma}_{n,r}^\flat(1)}T(\widetilde{\sigma})$ and
$\widetilde{\Gamma}_{n,r}^\flat(2)=\bigcup_{\sigma\in\Lambda_{n,r}^\flat(2)}S(\sigma)$. Write
\[
B_{1,r}:=\sum_{\sigma\in\Lambda_{n,r}^\flat(1)}\mathcal{E}_r(\sigma)^{t_r};\;B_{2,r}:=\sum_{\sigma\in\Lambda_{n,r}^\flat(2)}\mathcal{E}_r(\sigma)^{t_r}.
\]
In the following, we estimate $B_{1,r}$ and $B_{2,r}$ separately.

Note that $\{F_\sigma:\sigma\in\Lambda_{n,r}\}$ is a cover of the carpet $E$. For every $\widetilde{\sigma}\in\widetilde{\Gamma}_{n,r}^\flat(1)$, we have
$\sum_{\sigma\in T(\widetilde{\sigma})}\mu(F_\sigma)=\mu(F_{\widetilde{\sigma}})$. Further, by Remark \ref{r1}, we have, $|\sigma_L|-|\widetilde{\sigma}_L|\leq M_r$ for
every $\sigma\in T(\widetilde{\sigma})$. Note that $t_r\in (0,1)$. We deduce
\begin{eqnarray*}
\sum_{\sigma\in T(\widetilde{\sigma})}\mathcal{E}_r(\sigma)^{t_r}=\sum_{\sigma\in T(\widetilde{\sigma})}(\mu(F_\sigma)a_{\sigma_L^r})^{t_r}\geq\bigg(\sum_{\sigma\in T(\widetilde{\sigma})}\mu(F_\sigma)a_{\sigma_L^r}\bigg)^{t_r}
\geq \underline{a}^{rM_rt_r}\mathcal{E}_r(\widetilde{\sigma})^{t_r}.
\end{eqnarray*}
Let $H_{3,r}:=\underline{a}^{rM_rt_r}$. It follows that
\begin{eqnarray}\label{s2}
B_{1,r}=\sum_{\widetilde{\sigma}\in\widetilde{\Gamma}_{n,r}^\flat(1)}\sum_{\sigma\in T(\widetilde{\sigma})}\mathcal{E}_r(\sigma)^{t_r}
\geq H_{3,r}\sum_{\widetilde{\sigma}\in\widetilde{\Gamma}_{n,r}^\flat(1)}\mathcal{E}_r(\widetilde{\sigma})^{t_r}.
\end{eqnarray}
For every $\sigma\in\Lambda_{n,r}^\flat(2)$, we have, $S(\sigma)\subseteq G_1(\sigma)$. This and Lemma \ref{k3} yield
\begin{eqnarray}\label{s3}
B_{2,r}\geq \sum_{\sigma\in\Lambda_{n,r}^\flat(2)}H_{2,r}^{-1}\sum_{\widetilde{\sigma}\in S(\sigma)}\mathcal{E}_r(\widetilde{\sigma})^{t_r}=
H_{2,r}^{-1}\sum_{\widetilde{\sigma}\in\widetilde{\Gamma}_{n,r}^\flat(2)}\mathcal{E}_r(\widetilde{\sigma})^{t_r}.
\end{eqnarray}
Let $H_{4,r}:=\min(H_{2,r}^{-1},H_{3,r})$. Combining (\ref{s1})-(\ref{s3}), we obtain
\begin{eqnarray*}
\sum_{\sigma\in\Lambda_{n,r}}\mathcal{E}_r(\sigma)^{t_r}\geq
B_{1,r}+B_{2,r}\geq H_{4,r}\sum_{\widetilde{\sigma}\in\widetilde{\Gamma}_{n,r}^\flat}\mathcal{E}_r(\widetilde{\sigma})^{t_r}\gtrsim 1.
\end{eqnarray*}
This completes the proof of the lemma.
\end{proof}

\emph{Proof of Theorem \ref{mthm}} This is an easy consequence of Remark \ref{cardlambdajr}, Lemmas \ref{characterization}, \ref{k1}, \ref{k2} and \cite[Lemma 3.4]{Zhu:23}.

\subsection{Proof of Proposition \ref{mthm2}}
For every $n\geq 1$, we define
\begin{equation}\label{phits}
\Upsilon_n(t,s):=\sum_{\sigma\in\Phi_n}(p_{\sigma_L}q_{\sigma_R})^ta_{\sigma_L}^s,\;t\geq0,\;s\in\mathbb{R}.
\end{equation}
It is clear that one can replace $\Phi_n$ in (\ref{phits}) with $\Psi_n$. With some minor modifications of the proof of Lemmas \ref{lem1b}, \ref{lem1a}, one can obtain that
\[
\Upsilon_{p+k}(t,s)\asymp\Upsilon_p(t,s)\Upsilon_k(t,s).
\]
This allows us to define $\Upsilon(t,s):=\lim_{n\to\infty}\frac{1}{n}\log\Upsilon_n(t,s),\;t\geq0,s\in\mathbb{R}$. One can easily check that, for every $t\geq 0$ and $s\in\mathbb{R}$, $\Upsilon(t,s)$ is finite.
\begin{lemma}\label{zsg2}
For $t\in[0,\infty)$, there exists a unique number $\beta(t)$, such that $\Upsilon(t,\beta(t))=0$.
\end{lemma}
\begin{proof}
As in the proof of \cite[Lemma 5.2]{Fal:97}, for every $\epsilon>0$, we have
\begin{eqnarray*}
\epsilon\log\underline{a}\leq\Upsilon(t,s+\epsilon)-\Upsilon(t,s)\leq\epsilon\log\overline{a}.
\end{eqnarray*}
Thus, $\Upsilon(t,s)$ is strictly decreasing and continuous in $s$. Letting $\epsilon\to+\infty$, we obtain that $\lim\limits_{s\to\infty}\Upsilon(t,s)=-\infty$. Similarly, $\lim\limits_{s\to-\infty}\Upsilon(t,s)=+\infty$. Thus, the lemma follows by the continuity of $\Upsilon$ in $s$.
\end{proof}

\emph{Proof of Proposition \ref{mthm2}} Let $q\in[0,\infty)$ be given. Let $\beta(q)$ be as defined in Lemma \ref{zsg2}. For $\omega\in G^l$, as in Remark \ref{r2} (r2), let
\[
\Omega(\omega):=\{\tau\in G_y^*:\omega\times\tau\in\Phi_l\}.
\]
We have $b_{\omega_y} b_{\tau^-}\geq a_\omega>b_{\omega_y} b_\tau$. It follows that $b_\tau\asymp a_\omega/b_{\omega_y}$. Let $T_y(q)$ denote the $L^q$-spectrum for the projection of $\mu$ onto the y-axis: $\sum_{j\in G_y}q_j^qb_j^{T_y(q)}=1$. By induction, we have, $\sum_{\tau\in\Omega(\omega)}q_\tau^qb_\tau^{T_y(q)}=1$ for every $\omega\in G^*$. We deduce
\begin{eqnarray*}
\sum_{\sigma\in\Phi_l}(p_{\sigma_L}q_{\sigma_R})^qa_{\sigma_L}^{\beta(q)}&=&\sum_{\omega\in G^l}\sum_{\tau\in\Omega(\omega)}(p_{\omega}q_\tau)^qa_\omega^{\beta(q)}
\\&\asymp&\sum_{\omega\in G^l}p_{\omega}^qa_\omega^{\beta(q)}\sum_{\tau\in\Omega(\omega)}q_\tau^q(b_\tau^{T_y(q)}b_{\omega_y}^{T_y(q)}a_\omega^{-T_y(q)})\\&=&
\sum_{\omega\in G^l}p_{\omega}^qa_\omega^{\beta(q)-T_y(q)}b_{\omega_y}^{T_y(q)}\sum_{\tau\in\Omega(\omega)}q_\tau^qb_\tau^{T_y(q)}
\\&=&\sum_{\omega\in G^l}p_{\omega}^qa_\omega^{\beta(q)-T_y(q)}b_{\omega_y}^{T_y(q)}.
\end{eqnarray*}
This, along with (\ref{phits}), yields that
\[
\lim_{l\to\infty}\frac{1}{l}\log\sum_{\omega\in G^l}p_{\omega}^qa_\omega^{\beta(q)-T_y(q)}b_{\omega_y}^{T_y(q)}=\Upsilon(q,\beta(q))=0,
\]
which is equivalent to Feng-Wang's formula \cite[Theorem 2]{Feng:05}:
\[
\sum_{(i,j)\in G}p_{ij}^qa_{ij}^{\beta(q)-T_y(q)}b_j^{T_y(q)}=1.
\]
It follows that $\beta(q)=T(q)$ and the proof of the proposition is complete.



\begin{thebibliography}{35}

\bibitem{Bed:84} T. Bedford, Crinkly curves, Markov partitions and box dimensions in self-similar sets, PhD Thesis, University of Warwick, 1984.
\bibitem{Billings:84}P. Billingsley, Convergence of probability measures. John Wiley \& Sons, Inc, 1999.

\bibitem{Fal:97}K. J. Falconer, Techniques in fractal geometry, John Wiley \& Sons, Ltd., Chichester, 1997.
\bibitem{Feng:05}D.-J. Feng and Y. Wang, A class of self-affine sets and self-affine measures.  J. Fourier Anal. Appl.  11 (2005), 107-124.
\bibitem{GL:00} S. Graf and H. Luschgy, Foundations of quantization for probability distributions.  Lecture Notes in Math. vol. 1730, Springer, 2000.

\bibitem{GL:02}S. Graf and H. Luschgy, The quantization dimension of self-similar probabilities.  Math. Nachr.  241 (2002), 103-109.

\bibitem{GL:04}S. Graf and H. Luschgy, Quantization for probabilitiy measures with respect to the geometric mean error. Math. Proc. Camb. Phil. Soc. 136 (2004), 687-717.

\bibitem{GL:05}S. Graf and H. Luschgy, The point density measure in the quantization of self-similar probabilities. Math. Proc. Camb. Phil. Soc. 138 (2005), 513-31.
\bibitem{GL:12}S. Graf, H. Luschgy and  G. Pag\`{e}s, The local quantization behavior of absolutely continuous probabilities. Ann.  Probab. 40 (2012), 1795-1828.

\bibitem{GN:98} R. Gray and D. Neuhoff, Quantization. IEEE Trans. Inform. Theory 44 (1998), 2325-2383.


\bibitem{Hut:81} J. E. Hutchinson, Fractals and self-similarity. Indiana Univ. Math. J. 30 (1981), 713-47.
 \bibitem{KZ:15a}Marc Kesseb\"{o}hmer and Sanguo Zhu, Some recent developments in quantization of fractal measures, Fractal geometry and stochastics V, Progr. Probab., vol. 70, Birkh\"{a}user/Springer, Cham, 2015, pp. 105–120.

\bibitem{KZ:15} M. Kesseb\"{o}hmer and S. Zhu, On the quantization for self-affine measures on Bedford-McMullen carpets. Math. Z. 283 (2016), 39-58.
\bibitem{KNZ:22}M. Kesseb\"{o}hmer A. Niemann, and S. Zhu, Quantization dimensions of compactly supported probability measures via R\'{e}nyi dimensions. Trans. Amer. Math. Soc. 376 (2023), 4661-4678.
\bibitem{King:95}  J. F. King, The singularity spectrum for general Sierpi\'{n}ski carpets. Adv. Math. 116 (1995), 1-11.

\bibitem{Ko:22}I. Kolossv\'{a}ry, The $L^q$-spectrum of self-affine measures on sponges. J. London. Math. Soc.108 (2023), 666-701.


\bibitem{LG:92} S. P. Lalley and D. Gatzouras, Hausdorff and box dimensions of certain self-affine fractals. Indiana Univ. Math. J. 41 (1992), 533-568.


\bibitem{LM:02}L. J. Lindsay and R. D. Mauldin, Quantization dimension for conformal iterated function systems. Nonlinearity 15 (2002), 189-199.

\bibitem{Mcmullen:84} C. McMullen, The Hausdorff dimension of general Sierpi\'{n}ski carpets. Nagoya Math. J. 96 (1984), 1-9.

\bibitem{Olsen:98}L. Olsen, Self-affine multifractal Sierpi\'{n}ski sponges in $\mathbb{R}^d$, Pacific J. Math. 183 (1998), no. 1, 143–199.
\bibitem{Olsen:11}L. Olsen, Random self-affine multifractal Sierpi\'{n}ski sponges in $\mathbb{R}^d$, Monatsh. Math. 162 (2011), 89–117.
\bibitem{Peres:94b} Y. Peres, The self-affine carpets of McMullen and Bedford have infinite Hausdorff measure, Math. Proc. Camb. Phil. Soc. 116 (1994), 513-26.
\bibitem{Peres:00}Y. Peres and B. Solomyak, Existence of $L^q$ dimensions and entropy dimension for self-conformal measures, Indiana Univ. Math. J. 49 (2000), 1603-1621.
\bibitem{PK:01} K. P\"{o}tzelberger, The quantization dimension of distributions. Math. Proc. Camb. Phil. Soc. 131 (2001), 507-519.
\bibitem{PK:04}K. P\"{o}tzelberger, The quantization error of self-similar distributions.
Math. Proc. Camb. Phil. Soc. 137 (2004), 725-740.
\bibitem{Zhu:18}S. Zhu, Asymptotic order of the quantization error for a class of self-affine measures. Proc. Amer. Math. Soc. 146 (2018), 637-651.
\bibitem{Zhu:23}S. Zhu, Asymptotics of the quantization errors for some Markov-type measures with complete overlaps. J. Math. Anal. Appl. 528 (2023), 127585.
\bibitem{Zhu:25}S. Zhu, Asymptotic order of the quantization error for a class of self-similar measures with overlaps. Proc. Amer. Math. Soc. 153 (2025), 2115-2125.

\end{thebibliography}
\end{document}